\documentclass[10pt,a4paper]{article}
\usepackage{amstext}
\usepackage{array}
\usepackage{amsfonts}
\usepackage{amssymb, bm}
\usepackage{graphicx}
\usepackage{epstopdf}
\usepackage{algorithm}
\usepackage{algpseudocode}
\usepackage{url}
\usepackage{caption}
\usepackage{subfig}
\usepackage{csquotes}
\usepackage{tikz}
\usepackage{siunitx}
\usetikzlibrary{arrows}
\usetikzlibrary{decorations.markings}
\usetikzlibrary{shapes.geometric}
\usetikzlibrary{shapes.arrows}
\usetikzlibrary{fit,calc}
\usetikzlibrary{quotes,angles}
\usetikzlibrary{positioning, arrows.meta}
\usepackage{pgf}
\usepackage{pgfplots}
\usepackage{mathtools}
\usepackage{amsthm}
\usepackage{hyperref}
\usepackage{stmaryrd}
\usepackage{framed}
\usepackage{enumitem}
\usepackage[normalem]{ulem}
\usepackage{booktabs}
\usepackage{longtable}
\usepackage[dvipsnames]{xcolor}
\usepackage[capitalize,nameinlink]{cleveref}
\crefname{section}   {Section}   {Sections}
\crefname{subsection}{Subsection}{Subsections}
\Crefname{section}   {Section}   {Sections}
\Crefname{subsection}{Subsection}{Subsections}
\Crefname{figure}    {Figure}    {Figures}
\crefformat{equation}{\textup{#2(#1)#3}}
\crefrangeformat{equation}{\textup{#3(#1)#4--#5(#2)#6}}
\crefmultiformat{equation}{\textup{#2(#1)#3}}{ and \textup{#2(#1)#3}}
{, \textup{#2(#1)#3}}{, and \textup{#2(#1)#3}}
\crefrangemultiformat{equation}{\textup{#3(#1)#4--#5(#2)#6}}%
{ and \textup{#3(#1)#4--#5(#2)#6}}{, \textup{#3(#1)#4--#5(#2)#6}}{, and \textup{#3(#1)#4--#5(#2)#6}}
\Crefformat{equation}{#2Equation~\textup{(#1)}#3}
\Crefrangeformat{equation}{Equations~\textup{#3(#1)#4--#5(#2)#6}}
\Crefmultiformat{equation}{Equations~\textup{#2(#1)#3}}{ and \textup{#2(#1)#3}}
{, \textup{#2(#1)#3}}{, and \textup{#2(#1)#3}}
\Crefrangemultiformat{equation}{Equations~\textup{#3(#1)#4--#5(#2)#6}}%
{ and \textup{#3(#1)#4--#5(#2)#6}}{, \textup{#3(#1)#4--#5(#2)#6}}{, and \textup{#3(#1)#4--#5(#2)#6}}

\crefdefaultlabelformat{#2\textup{#1}#3}
\crefname{proposition}{Proposition}{Propositions}
\Crefname{proposition}{Proposition}{Propositions}
\crefname{definition} {Definition} {Definitions}
\Crefname{definition} {Definition} {Definitions}
\crefname{theorem}    {Theorem}    {Theorems}
\Crefname{theorem}    {Theorem}    {Theorems}
\crefname{remark}     {Remark}     {Remarks}
\Crefname{remark}     {Remark}     {Remarks}
\crefname{assumption} {Assumption} {Assumptions}
\Crefname{assumption} {Assumption} {Assumptions}

\newcommand{\vertiii}[1]{{\left\vert\kern-0.25ex\left\vert\kern-0.25ex\left\vert #1
   \right\vert\kern-0.25ex\right\vert\kern-0.25ex\right\vert}}

\newcommand{\bs}[1]{\boldsymbol{#1}}
\newcolumntype{T}[1]{>{\ttfamily}p{#1}}
\newtheorem{theorem}{Theorem}[section]

\newtheorem{proposition}[theorem]{Proposition}

\theoremstyle{definition}
\newtheorem{remark}[theorem]{Remark}

\usepackage{geometry}
\date{ \today}
\title{
Model order reduction of HDG discretizations for harmonic wave propagation
}
\author{H\'el\`ene Barucq\thanks{Project-Team Makutu, Centre Inria de l'Université de Bordeaux, Universit\'e de Pau et des Pays de l'Adour, LMAP, CNRS UMR 5142, Pau, France.\\ Emails: \texttt{helene.barucq@inria.fr}, \texttt{florian.faucher@inria.fr}, \texttt{ha.howard@inria.fr}, \texttt{ishak.tifouti@inria.fr}.}, \,
Florian Faucher\footnotemark[1], \, 
Ha Pham\footnotemark[1], \, 
Ishak Tifouti\footnotemark[1]}

\begin{document}
\maketitle

\begin{abstract}
We propose and analyze a projection-based model reduction method for the rapid and reliable solution of time-harmonic acoustic wave propagation problems over a frequency band, with particular emphasis on seismic applications. The high-fidelity model is based on a hybridizable discontinuous Galerkin (HDG) discretization, formulated as a generalized saddle-point system coupling the volume and trace unknowns.  The basis construction combines residual-based greedy sampling of the frequency parameter, enrichment with frequency derivatives of the HDG solution, and proper orthogonal decomposition (POD) compression. We establish the well posedness of both the high-fidelity and reduced problems and characterize their stability with inf-sup constants. We investigate the interaction between static condensation, inherent to the HDG discretization, and projection onto the reduced spaces, characterizing its implications for the reduced formulation and approximation errors. Since these operations do not commute in general, we derive a priori error estimates for both the trace and volume variables to quantify the gap between the two operations. We further identify a lifting condition on the volume reduced space that eliminates these gaps and yields a trace error estimate independent of the volume approximation error. Numerical experiments on the two-dimensional Marmousi model and a three-dimensional heterogeneous medium, with more than $10^{7}$ unknowns, illustrate the theoretical estimates and show that the reduced model achieves speedups of several orders of magnitude compared to the high-fidelity solver.
\end{abstract}

\section{Introduction}
\label{sec:intro}

In the context of time-harmonic wave propagation, a frequency sweep is used to describe the dependence of the wavefield on the frequency and therefore requires the repeated solution of the forward problem for multiple parameter values. For large-scale simulations, particularly in three dimensions, the computational cost of such multi-query evaluations can rapidly become prohibitive. Hybridizable discontinuous Galerkin (HDG) methods \cite{cockburn2009unified} provide an attractive discretization framework for such problems, combining the high-order accuracy and local approximation properties of discontinuous Galerkin methods with a reduced global coupling structure. In particular, HDG method distinguishes between element-wise volume variables and a hybrid trace variable defined on the mesh skeleton, which couples neighboring elements: the volume unknowns are eliminated cell by cell, so that the global system involves the trace variable only. The approach has been sucessefully applied to a large class of problems, for instance in continuum mechanics \cite{nguyen2012hybridizable} and to wave propagation \cite{nguyen2011high,fernandez2018hybridized}, where it yields optimal convergence rates for all approximated fields. For the Helmholtz equation, error analyses are available in \cite{griesmaier2011error,cui2014analysis}. The structure of HDG methods is well suited to large-scale and parallel computations; however, the repeated solution of the resulting high-dimensional systems over a range of frequencies can still be computationally demanding. Model order reduction (MOR) for parametric problems provides a means of alleviating this cost by reducing the marginal cost of repeated solves, thereby enabling denser sampling of the frequency range and the exploration of broader frequency bands within an acceptable computational budget.

\subsection{Relation to previous works}
The MOR approach considered in this work is projection-based, also referred to as intrusive: the discretized full-order model (FOM) is projected onto a low-dimensional approximation space, resulting in a reduced order model (ROM) whose online evaluation is much cheaper than that of the full system. The construction of the reduced space depends on the structure of the underlying problem, the available computational resources, and the objectives of the reduction. Within the parametrized partial differential equation (pPDE) community, the variational formulation of the PDE provides an attractive framework for analyzing the approximation and stability properties of the resulting reduced model; see \cite{quarteroni2014reduced,benner2017model,quarteroni2011certified} for a comprehensive treatment of reduced basis methods for parametrized PDEs. Reduced spaces are commonly constructed from solution snapshots using Proper Orthogonal Decomposition (POD) \cite{bakewell1967viscous,sirovich1987turbulence}, which provides a basis that is optimal in the mean-square sense over the prescribed set of sampled solutions. Alternatively, weak-greedy algorithms construct the reduced space adaptively by selecting parameter values sequentially, with each new sample driven by a computationally inexpensive a posteriori error estimator or residual-based indicator \cite{prud2002reliable}.

In contrast, the model reduction literature for dynamical, linear time-invariant, and algebraic systems often constructs reduced spaces directly from the algebraic structure of the discrete operators, for instance through Krylov subspaces or moment-matching techniques, rather than through the explicit sampling of solution fields. For instance, Asymptotic Waveform Evaluation (AWE) \cite{slone2003broadband} and its well-conditioned variant, WCAWE \cite{slone2003well}, construct reduced models by matching solution derivatives, thereby ensuring rational interpolation up to a prescribed order around selected expansion points. The multivariate extension MWCAWE further enables joint moment matching with respect to several parameters \cite{rumpler2023mwcawe}. Closely related approaches include rational or Padé-type approximations of solution maps \cite{bonizzoni2020least}, as well as shift-adaptive rational Krylov methods for wave propagation \cite{druskin2010adaptive}. In contrast to POD- and greedy-based constructions, which typically rely on global sampling of the parameter domain, these methods enrich the reduced space through explicitly computed moments at a comparatively small number of expansion points. Their effectiveness consequently depends on the efficient computation of these moments, which relies on the linearity of the underlying model and the regularity of the associated parameter-to-solution map. We refer to \cite{benner2015survey,benner2017model} for a survey of projection-based model reduction methods for parameterized systems.

More generally, several strategies have been proposed to enrich reduced spaces with information tailored to a specific objective. Derivative-based enrichment can improve local approximation properties around expansion points; see, for example, \cite{Haasdonk}. Regularization-aware enrichment has been investigated in reduced order optimization \cite{kartmann2024adaptive}, while adjoint-based enrichment can improve the accuracy of reduced models for prescribed quantities of interest through goal-oriented model reduction \cite{yano2021globally,zahr2015progressive,tifouti2025construction}. Finally, the approaches discussed above are referred to as linear model reduction techniques in the sense that the approximation is sought within a single linear subspace of fixed dimension. This assumption may become restrictive when the solution set can not be accurately represented by a single low-dimensional linear space. One possible alternative is to construct a collection of local reduced spaces associated with a partition of the parameter domain \cite{eftang2010hp,dihlmann2011model,Haasdonk,tifouti2025registration}. Alternatively, genuinely nonlinear approximation strategies can be employed, for instance through nonlinear manifold approximations based on autoencoders \cite{lee2020model,barnett2023neural,pichi2023artificial}. In the present time-harmonic formulation of wave propagation, however, as explained in \cite{henriquez2026fast}, the favorable approximation properties associated with the frequency-domain solution map make such nonlinear constructions unnecessary for the class of problems considered here. We therefore restrict our methodology to linear reduced spaces throughout this work.

MOR has also been successfully applied to seismic and wave propagation problems. In particular, model reduction techniques have been applied to frequency-domain and normal-mode formulations in wave propagation and seismic imaging \cite{hawkins2023model,hawkins2025model}. Reduced models have also been employed as computational surrogates within inversion procedures \cite{besset2025model,bui2008model,borcea2014model}. Extensions to saddle-point and mixed formulations, which are particularly relevant to the HDG discretization considered here, are discussed in \cite{gerner2012certified,gernerVeroy2012certified}. The  latter is relevant to the present work, as the HDG formulation considered in this work can be expressed as a generalized saddle-point problem, involving distinct volume and hybrid variables together with their coupling through the mesh skeleton.

Our point of departure is the HDG formulation introduced in \cite{faucher2020adjoint} for an acoustic damped wave-propagation problem. The introduction of a trace variable on the mesh skeleton enables the element-wise elimination of the volume unknowns and leads, after static condensation, to a global problem involving only the hybrid unknown. This structure reduces the size of the globally coupled system. From the perspective of model reduction, however, the hybrid structure of HDG raises specific issues. The formulation involves distinct volume and hybrid variables coupled through the mesh interfaces, and the choice between reducing the condensed and non-condensed formulations has important consequences for the resulting reduced model. In particular, the condensed system involves the inverse of parameter-dependent matrices, which destroys the affine parameter dependence required for an efficient offline/online decomposition. To address this issue, the authors of \cite{vidal2015model} construct a reduced problem for an HDG system while preserving an affine structure by eliminating only the flux variables. The effectiveness of their strategy is demonstrated for the stochastic evaluation of a quantity of interest associated with a Helmholtz problem. 

Another related class of model reduction methods consists of localized model reduction approaches, in which local reduced spaces are constructed on subdomains or components and subsequently coupled through suitable interface conditions. Examples include the reduced basis element method, which combines reduced basis approximation with domain decomposition \cite{maday2002reduced}, and the static condensation reduced basis element method, which combines reduced basis method with component mode synthesis \cite{phuong2013static}.
The interface space can be reduced as well, either empirically \cite{eftang2013port} or from the dominant eigenfunctions of a transfer operator \cite{smetana2016optimal}; the latter construction is related to the generalized finite element method \cite{babuska2011optimal}.  We refer to \cite{buhr2020localized} for a survey of localized model reduction methods.

\subsection{Contributions}
In our work, we build the reduced basis through weak-greedy sampling, followed by POD compression, and it is subsequently enriched with Fréchet derivatives of the parametrized HDG solution map. Including derivatives provides local information on the variation of the solution map and is expected to improve its representation in the vicinity of the corresponding expansion points. We therefore combine global adaptive sampling with local derivative information, with the aim of improving the efficiency of frequency-domain reduction without increasing the number of full-order training solves. 

Similarly to \cite{vidal2015model}, we construct the reduced order model by projecting the non-condensed HDG formulation, thereby preserving the affine parameter dependence of the underlying discretized problem. Our approach, however, differs in several aspects. In contrast to \cite{vidal2015model}, we retain the complete HDG formulation, including the flux variables, which are essential for the inversion applications targeted in this work. Furthermore, we construct separate reduced spaces for the volume and hybrid trace variables. This distinction is central to our methodology and plays an important role in the well-posedness and error analyses developed below.

Beyond the construction of the reduced model, we provide a mathematical analysis of the well-posedness of both the full-order HDG formulation and its reduced counterpart. The use of separate reduced spaces for the volume and hybrid variables plays a central role in this analysis, as it (i) provides additional flexibility in the choice of reduced trial and test functions and (ii) allows the HDG structure to be exploited when deriving the corresponding energy identities. We further derive a priori error estimates for the reduced approximation that explicitly quantify the propagation of the hybrid approximation error to the volume variables.

Finally, we show that the error estimates provide a means of constructing an enrichment strategy based on a lifting method. In particular, we show that incorporating  a lifting of the trace-coupling contribution into the volume reduced space allows the error associated with the hybrid variable to be completely decoupled from the approximation error associated with the volume reduced space. As a consequence, the lifting strategy prevents the approximation of the hybrid variable from being degraded by the choice of the volume reduced basis and therefore allows it to inherit the same approximation properties as the corresponding full-order HDG approximation.

As a consequence, the resulting framework combines the computational advantages of HDG with a reduction strategy specifically adapted to its hybrid structure. We assess the performance of our methodology on two- and three-dimensional heterogeneous acoustic wave propagation model problems.  This work constitutes a baseline for the incorporation of reduced order models into full waveform inversion, and we therefore restrict the parametrization to the frequency for a fixed medium; the treatment of medium variability is left for future work. The individual ingredients of the proposed methodology are well established in the reduced basis literature. The novelty of our work lies in their combination within the HDG setting, together with the analysis we provide, which establishes the well-posedness of both the full-order and reduced problems and shows how key structural properties of the HDG scheme are inherited at the reduced level through our methodology. 

The remainder of the paper is organized as follows. Section~\ref{sec:model} presents the time-harmonic acoustic model and the treatment of the boundary conditions. Section~\ref{sec:discretization} introduces the HDG discretization, discusses its condensed and uncondensed forms, and establishes the well-posedness of the resulting systems together with the analyticity of the solution map with respect to the frequency. Section~\ref{sec:ROM} describes the construction of the reduced bases, the Galerkin projection, and the associated error analysis, including estimates that quantify the gap between projecting the condensed system and condensing the projected one, illustrating that the two operations do not commute. In~\cref{sec:marmousi,sec:3d}, we conduct extensive numerical investigations of the different ingredients of our methodology in two- and three-dimensional heterogeneous media, and Section~\ref{sec:conclusion} concludes the paper.


\section{Acoustic wave propagation}
\label{sec:model}
We consider the propagation of damped time-harmonic waves in acoustic media 
in the domain $\Omega \subset \mathbb{R}^d$, $d=2,3$, with Lipschitz boundary 
$\Gamma = \Gamma_N \cup \Gamma_A$, $\Gamma_N\cap\Gamma_A=\emptyset$. 
The medium is characterized by its density $\rho\in L^\infty(\Omega)$ and
bulk modulus $\kappa\in L^\infty(\Omega)$, assumed piecewise constant and satisfying $0<\rho_{\min}\le\rho\le\rho_{\max}$ and 
$0<\kappa_{\min}\le\kappa\le\kappa_{\max}$. 
This ensures that $\rho$ and $\kappa$ have well-defined, 
bounded traces on $\Gamma$, as required by the boundary condition below.
We introduce the complex frequency $\sigma =  i\omega - s \in \mathbb{C}$,
where $\omega\in\mathbb{R}^*_+$ is the angular frequency and $s>0$ 
acting as a damping parameter, e.g., \cite{Faucher2023viscoacoustic}. 
The first-order velocity--pressure formulation of the time-harmonic
acoustic system reads:\\
\noindent
Given the scalar source $f\in L^2(\Omega)$, find $p \in H^1(\Omega ; \mathbb{C})$ and
$\boldsymbol{v} \in H_{\operatorname{div}}( \Omega ; \mathbb{C}^d)$ such that
\begin{equation}
\left\{
\begin{aligned}
- \dfrac{\sigma}{\kappa(\boldsymbol{x})}\,p(\boldsymbol{x})
   + \nabla\cdot\boldsymbol{v}(\boldsymbol{x}) &= f(\boldsymbol{x}),
   &&\text{in }\Omega,\\
- \sigma\,\rho(\boldsymbol{x})\,\boldsymbol{v}(\boldsymbol{x})
   + \nabla p(\boldsymbol{x}) &= \boldsymbol{0},
   &&\text{in }\Omega,\\
\boldsymbol{v}(\boldsymbol{x})\cdot\boldsymbol{\nu} &= 0,
   &&\text{on }\Gamma_N,\\
-\bigl(c(\boldsymbol{x})\,\rho(\boldsymbol{x})\bigr)^{-1} p(\boldsymbol{x})
   + \boldsymbol{v}(\boldsymbol{x})\cdot\boldsymbol{\nu} &= 0,
   &&\text{on }\Gamma_A,
\end{aligned}
\right.
\label{eq:acoustic_strong}
\end{equation}
where $c(\boldsymbol{x}) = \sqrt{\kappa(\boldsymbol{x})/\rho(\boldsymbol{x})}$ 
is the wave speed. \\
We note that, since $\rho$ and $\kappa$ are assumed piecewise constant and bounded,
$(c\rho)^{-1}|_\Gamma=1/\sqrt{\kappa\rho}\,|_\Gamma\in L^\infty(\Gamma)$; together with the traces
$p|_\Gamma\in H^{1/2}(\Gamma)$ and $(\boldsymbol{v}\cdot\boldsymbol{\nu})|_\Gamma\in H^{-1/2}(\Gamma)$,
this makes the boundary conditions well defined in $H^{-1/2}(\Gamma)$.

System~\cref{eq:acoustic_strong} is derived from the first-order time-domain acoustic equations by seeking time-harmonic solutions with complex frequency $\sigma=i\omega-s$, see~\cite[Appendix~A]{Faucher2023viscoacoustic}. Indeed, the corresponding temporal dependence is $e^{\sigma t}=e^{-st}e^{i\omega t}$, which combines oscillations at angular frequency $\omega$ with exponential damping at rate $s$. The resulting frequency-domain formulation is  referred to as the Laplace--Fourier domain~\cite{Shin2009}.
The Neumann condition on $\Gamma_N$ models a reflecting boundary. 
When $\Omega$ represents a portion of an unbounded medium, artificial absorbing 
boundary conditions are prescribed on $\Gamma_A$: following 
\cite{Engquist1977,faucher2020adjoint}, we employ the first-order absorbing 
boundary condition (ABC) in the last line of~\cref{eq:acoustic_strong}.
A key advantage of this first-order formulation is that the ABC is
expressed directly in terms of the physical fields $(p,\boldsymbol{v})$ and is
independent of the frequency $\sigma$. 
This property is particularly attractive in the reduced basis framework where 
reduced operators must be assembled based on their dependence on the frequency.

We further define the reciprocal quantity $\mu$ 
and the notation $\mathbf{m}$ that will be used 
in the rest of the manuscript such that, 
\begin{equation}
   \mu(\bs{x}) := \frac{1}{\kappa(\bs{x})},
   \qquad
   \mathbf{m} := (\rho,\mu).
\label{eq:media_def}
\end{equation}

\section{HDG discretization}
\label{sec:discretization}

We briefly present the hybridizable discontinuous Galerkin discretization
of  \eqref{eq:acoustic_strong}, following \cite{faucher2020adjoint}. Let
$\mathcal{T}_h$ be a conforming mesh of $\Omega$ with cells $\{K_e\}$ and
skeleton $\Sigma = \bigcup_e \partial K_e$ composed of faces
$\mathfrak{f}$. The HDG method introduces a hybrid unknown $\lambda_h$ on
$\Sigma$ and seeks
$(p_h,\bs{v}_h,\lambda_h)\in \mathcal{V}_h := W_h \times\bs{W}_h\times M_h$, where 
\begin{equation}
   W_h := \bigl\{ w\in L^2(\Omega) : w|_{K_e}\in\mathbb{P}_p(K_e)\ \ \forall K_e\in\mathcal{T}_h \bigr\},
   \qquad
   \boldsymbol{W}_h := (W_h)^d,
   \qquad
   M_h := \bigl\{ q\in L^2(\Sigma) : q|_{\mathfrak{f}}\in\mathbb{P}_p(\mathfrak{f})\ \ \forall \mathfrak{f}\subset\Sigma \bigr\}.
   \label{eq:hdg_spaces}
\end{equation}
\begin{itemize}
    \item On each cell $K_e$, let $\{\phi_j\}_{j=1}^{N_{\mathrm{dof}}^{(e)}}$
    be a basis of $\mathbb{P}_{p}(K_e)$, used to define the
    broken pressure space $W_h$ and velocity space $\bs{W}_h$:
    \begin{equation}
        p_h|_{K_e} = \sum_{j} \mathrm{p}_j^{(e)}\,\phi_j,
        \qquad
        v_{x,h}|_{K_e} = \sum_{j} \mathrm{v}_{x,j}^{(e)}\,\phi_j,
        \qquad
        v_{y,h}|_{K_e} = \sum_{j} \mathrm{v}_{y,j}^{(e)}\,\phi_j,
    \end{equation}
    with local unknown vector
    $U_e := (\mathrm{p}^{(e)},\mathrm{v}_x^{(e)},\mathrm{v}_y^{(e)})^T$.

    \item On each face $\mathfrak{f}\subset\Sigma$, let
    $\{\chi_k\}_{k=1}^{\widehat{N}_{\mathrm{dof}}^{(\mathfrak{f})}}$ be a
    basis of $\mathbb{P}_{p}(\mathfrak{f})$,
    defining the trace space $M_h$, with
    $\lambda_h|_{\mathfrak{f}} = \sum_k \lambda_{h,k}^{(\mathfrak{f})}\,\chi_k$
    and global face vector $\Lambda$.
\end{itemize}

\noindent
The restriction operator $\mathcal{R}_e$ maps $\Lambda$ to its
restriction on $\partial K_e$: $\mathcal{R}_e\Lambda = \Lambda|_{\partial K_e}$.\\
We denote by
\begin{equation}
   N_U := (d+1)\sum_{e} N_{\mathrm{dof}}^{(e)},
   \qquad
   N_\Lambda := \sum_{\mathfrak{f}\subset\Sigma} \widehat{N}_{\mathrm{dof}}^{(\mathfrak{f})},
   \label{eq:dof_counts}
\end{equation}
the numbers of volume and trace degrees of freedom, which are the sizes of the
global vectors $U=(U_e)_e$ and $\Lambda$. 

\subsection{Local matrices}
\label{sec:local_matrices}

On each cell $K_e$, the local matrix reads
\begin{equation}
\mathbb{A}_e(\sigma, \rho, \mu, \tau)=
\begin{pmatrix}
-(\sigma\mu\,\phi_i,\phi_j)_{K_e}
+\tau(\phi_i,\phi_j)_{\partial K_e}
& (\partial_x\phi_i,\phi_j)_{K_e}
& (\partial_y\phi_i,\phi_j)_{K_e}\\[2pt]
-(\phi_i,\partial_x\phi_j)_{K_e}
&-(\sigma\rho\phi_i,\phi_j)_{K_e} & 0\\[2pt]
-(\phi_i,\partial_y\phi_j)_{K_e}
& 0
&-(\sigma\rho\phi_i,\phi_j)_{K_e}
\end{pmatrix},
\label{eq:Ae_def}
\end{equation}
where $\tau>0$ is the stabilization parameter defined in \cref{sec:wellposedness}. The coupling matrices are
\begin{equation}
\mathbb{C}_e(\tau)=
\begin{pmatrix}
-\tau(\chi_k,\phi_j)_{\mathfrak{f}}\\[2pt]
(\chi_k,\phi_j\nu_{x,\mathfrak{f}})_{\mathfrak{f}}\\[2pt]
(\chi_k,\phi_j\nu_{y,\mathfrak{f}})_{\mathfrak{f}}
\end{pmatrix}_{\mathfrak{f}=\mathfrak{f}_1,\ldots,\mathfrak{f}_{N_{\mathrm{face}}}},
\qquad
\mathbb{B}_e^{(\mathfrak{f})}(\tau)
=
\Bigl(
\tau(\chi_k,\phi_j)_{\mathfrak{f}}\;\;
(\chi_k,\phi_j\nu_{x,\mathfrak{f}})_{\mathfrak{f}}\;\;
(\chi_k,\phi_j\nu_{y,\mathfrak{f}})_{\mathfrak{f}}
\Bigr),
\label{eq:BCe_def}
\end{equation}
the local right-hand side is
$\mathbb{S}_e = \bigl((f,\phi_j)_{K_e},\,0,\,0\bigr)^T$, and the
face block $\mathbb{L}_e$ is block-diagonal with
\begin{equation}
\mathbb{L}_e^{(\mathfrak{f}^{\mathrm{I}})}(\tau)
= - \tau\,(\chi_k,\chi_j)_{\mathfrak{f}^{\mathrm{I}}},
\qquad
\mathbb{L}_e^{(\mathfrak{f}^{\mathrm{B}})} (\tau, \rho, \mu)
= - \bigl(c_e + \tau\bigr)(\chi_k,\chi_j)_{\mathfrak{f}^{\mathrm{B}}},
\qquad
c_e := \sqrt{\mu\,\rho^{-1}},
\label{eq:Le_def}
\end{equation}
where $\mathrm{I}$ and $\mathrm{B}$ denote interior and boundary faces.

\subsection{Condensed formulation}
\label{sec:condensed}
The classical formulation of the HDG system is based on the Schur complement which  yields the local system
\begin{equation}
    \mathbb{A}_e\,U_e + \mathbb{C}_e\,\mathcal{R}_e\Lambda = \mathbb{S}_e.
\label{eq:local}
\end{equation}
We express the volume variables using global surface variables:
\begin{equation}
    U_e = \mathbb{A}_e^{-1}
    \bigl(\mathbb{S}_e - \mathbb{C}_e\,\mathcal{R}_e\Lambda\bigr).
\label{eq:backsubstitution}
\end{equation}
Substituting into the global transmission condition on $\Sigma$ gives the
condensed system
\begin{equation}
    \mathbb{K}\,\Lambda = \mathbb{F}, 
    \label{eq:condensed}
\end{equation}
where 
\begin{equation}
    \mathbb{K} :=
    \sum_e \mathcal{R}_e^T
    \bigl(\mathbb{L}_e - \mathbb{B}_e\,\mathbb{A}_e^{-1}\mathbb{C}_e\bigr)
    \mathcal{R}_e,
    \qquad
    \mathbb{F} :=
    -\sum_e \mathcal{R}_e^T\,\mathbb{B}_e\,\mathbb{A}_e^{-1}\mathbb{S}_e.
\label{eq:K_def}
\end{equation}
This is the standard formulation for high-fidelity solves which makes the HDG system attractive since the size of the global system is
determined by face degrees of freedom and the local variables $U_e$ are obtained through parallel computing.\\
Note that the operator $\mathbb{K}$ is no longer affine with respect to the frequency $\sigma$ since it involves $\mathbb{A}_e^{-1}$. The affine structure is however useful in the model reduction framework as it allows for an efficient offline-online decomposition of the computations (Cf. \cref{sec:galerkin}). 
Next we present an equivalent formulation of HDG that preserves the affine structure with respect to $\sigma$.

\subsection{Uncondensed formulation}
\label{sec:uncondensed}
Note that the condensed formulation \cref{eq:condensed} is used to compute the high-fidelity
solutions, however the local elimination of the volume variables introduces the
inverse of $\mathbb{A}_e$ and consequently destroys the affine dependence of
the global operator $\mathbb{K}$ on the frequency. For the reduced order model, we
therefore keep the volume and trace variables as independent unknowns.
Together with the global transmission condition
$$
\sum_e \mathcal{R}_e^T
\bigl(\mathbb{B}_e U_e+\mathbb{L}_e\mathcal{R}_e\Lambda\bigr)=0,
$$
the local systems \eqref{eq:local} can be written in the equivalent
uncondensed block form
\begin{equation}
    \underbrace{\begin{pmatrix}
    \mathbb{A}& \mathbb{C}\\
    \mathbb{B} & \mathbb{L}
    \end{pmatrix}}_{\mathbb{G}}
    \begin{pmatrix} U \\ \Lambda \end{pmatrix}
    =
   \underbrace{\begin{pmatrix} \mathbb{S} \\ 0 \end{pmatrix}}_{\mathbb{R}_S},
\label{eq:full_block}
\end{equation}
where $\mathbb{A}=\mathrm{diag}_e(\mathbb{A}_e)$,
$\mathbb{C}=\mathrm{diag}_e(\mathbb{C}_e\mathcal{R}_e)$, and similarly for
$\mathbb{B}$, $\mathbb{L}$, and $\mathbb{S}$. In contrast to the condensed
operator \eqref{eq:K_def}, all blocks in \eqref{eq:full_block} preserve the
affine dependence on the frequency $\sigma$. 

\begin{remark}
The formulation considered here is related to the approach of
\cite{vidal2015model,vidal2016empirical}, where a two-field formulation involving only
$(p_h,\lambda_h)$ is used and the flux $\bs{v}_h$ can be recovered by
back-substitution. While their formulation preserves the affine frequency
dependence and reduces the dimension of the variables involved during the
offline stage; Recovering the flux variables when they are required adds
an additional online computation involving the high-fidelity volume
operators. In our setting, the velocity field is part of the quantities of
interest and is therefore retained explicitly in the reduced formulation.

Accordingly, our approach differs from \cite{vidal2015model,vidal2016empirical} in two respects. First, we retain the
three-field formulation $(p_h,\bs{v}_h,\lambda_h)$ and compensate for the
increased offline cost through parallel computation. Second, we construct
distinct reduced spaces for the volume variables
$u_h=(p_h,\bs{v}_h)$ and the trace variable $\lambda_h$. This separation is
required with the well-posedness analysis as explained in Remark
\ref{rem:wellposdeness_ROM}, which relies on independently constructed
approximation spaces for the volume and trace variables. It also preserves
the block structure of the high-fidelity formulation in the reduced system,
as discussed in \cref{sec:galerkin}.
\end{remark}

\begin{remark}
\label{rem:local_proj}
We note that the uncondensed system \eqref{eq:full_block} is neither assembled nor solved
globally in our implementation. Indeed, solving the system \eqref{eq:full_block} would
introduce unnecessary global volume degrees of freedom that the HDG method aims to mitigate.  During the training stage, the condensed
system \eqref{eq:condensed} is therefore solved to compute the high-fidelity snapshots.
During the projection stage, only the local contributions
$\mathbb{A}_e$, $\mathbb{B}_e$, $\mathbb{C}_e$, $\mathbb{L}_e$, and
$\mathbb{S}_e$ are accessed element-wise and projected onto the corresponding
reduced spaces. The global uncondensed operator therefore never needs to be
assembled. This element-wise construction preserves the HDG structure,
avoids unnecessary memory requirements, and enables the projection operations
to be performed in parallel.
\end{remark}

\subsection{Well-posedness}
\label{sec:wellposedness}
In this section we prove the well posedness of the HDG system \eqref{eq:full_block}. Authors in \cite{du2019invitation} prove well-posedness of the HDG system for
the undamped Helmholtz equation $s=0$, and \cite{barucq2023construction} for the convected
Helmholtz equation without damping, both works rely on a G\r{a}rding inequality together with a duality argument. In our work, we consider a damped Laplace-Fourier setting, the fact that $s>0$ makes the operators sign definite. This makes the proofs simpler and requires no a priori error estimates. 

First, we prove the well-posedness of the global system  \eqref{eq:full_block}, then we prove the invertibility of
$\mathbb{A}_e$, on which the condensation  \eqref{eq:backsubstitution}
relies. 

We first define the sesquilinear forms associated with the
systems  \eqref{eq:local}, \eqref{eq:full_block}
and  \eqref{eq:condensed}. In what follows we write $u_h:=(p_h,\boldsymbol v_h)$, and $U$, $\Lambda$  their corresponding
coefficient vectors. Similarly, for a generic function $x=(u,\lambda)\in\mathcal{V}_h$,
 $X=(U,\Lambda)$ stands for its vector coefficients.\\
We denote by $\mathbb{M}_U$ and $\mathbb{M}_\Lambda$ the norm matrices  \[(u,w)_{L^2(\Omega)} = W^{*}\mathbb{M}_U U, \qquad  (\lambda,q)_{L^2(\Sigma)} = Q^{*}\mathbb{M}_\Lambda\Lambda \] 
and $\mathbb{M}:=\mathrm{diag}(\mathbb{M}_U,\mathbb{M}_\Lambda)$.

On a cell $K_e$, for
$u_e=(p,\boldsymbol v)$ and $w_e=(\phi,\boldsymbol\psi)$ in
$W_h(K_e)\times\boldsymbol W_h(K_e)$, we set
\begin{equation}
  a_e(u_e,w_e) := W_e^{*}\,\mathbb{A}_e\,U_e ,
  \qquad
  a(u,w) := \sum_{e} a_e(u_e,w_e)
          = W^{*}\,\mathbb{A}\,U ,
\label{eq:form_ae}
\end{equation}
with $\mathbb{A}_e$ given by  \eqref{eq:Ae_def}. 

For the uncondensed system, let $x=(u,\lambda)$ and $y=(w,q)$ in
$\mathcal{V}_h=\bigl(W_h\times\boldsymbol W_h\bigr)\times M_h$, and set
\begin{equation}
  \mathcal{A}(x,y) := Y^{*}\,\mathbb{G}\,X 
\label{eq:form_calA}
\end{equation}
$\mathbb{G}$ being the block matrix of  \eqref{eq:full_block}.

Finally, for the condensed system we set
\begin{equation}
  \mathcal{K}(\lambda,q) := Q^{*}\,\mathbb{K}\,\Lambda ,
\label{eq:form_calK}
\end{equation}
with $\mathbb{K}$ given by  \eqref{eq:K_def}.

We equip $\mathcal{V}_h$ with the norm
\begin{equation}
  \|x\|_{\mathcal{V}_h}^2 := \|p\|_{L^2(\Omega)}^2 + \|\boldsymbol v\|_{L^2(\Omega)}^2 + \|\lambda\|_{L^2(\Sigma)}^2 = X^{*}\,\mathbb{M}\,X ,
\label{eq:norm}
\end{equation}
and introduce the  operator
\begin{equation}
  T : \mathcal{V}_h\to \mathcal{V}_h ,
  \qquad
  T(p,\boldsymbol v,\lambda) := (p,\ \boldsymbol v,\ -\lambda) ,
\label{eq:Tdef}
\end{equation}
that satisfies $T^2=I$.

\begin{proposition}[T-coercivity of $\mathcal{A}$]
    Assume $s>0$, $\tau\ge\tau_{\min}>0$ and
$\mu\ge\mu_{\min}>0$. Then, $\exists \beta >0 $ such that for every $x\in \mathcal{V}_h$ we have, 
\begin{equation}
  \operatorname{Re}\,\mathcal{A}(x,Tx)\ \ge\ \beta\,\|x\|_{\mathcal{V}_h}^2 .
\label{eq:Tcoerc}
\end{equation}
\end{proposition}
 
\begin{proof}
Let $x=(p,\boldsymbol v,\lambda), y=(\phi,\boldsymbol\psi,q) \in \mathcal{V}_h$ , we have
\begin{equation}
\begin{aligned}
  \mathcal{A}(x,y) = \;&
  \sum_{K_e}\Bigl(-\sigma\,\rho(\boldsymbol v,\boldsymbol\psi)_{K_e}
    -(p,\nabla\cdot\boldsymbol\psi)_{K_e}
    +\langle\lambda,\boldsymbol\psi\cdot\boldsymbol \nu\rangle_{\partial K_e}\Bigr)
  \\
  &+\sum_{K_e}\Bigl(-\sigma\,\mu\,(p,\phi)_{K_e}
    +(\nabla\cdot\boldsymbol v,\phi)_{K_e}
    +\tau\,\langle p-\lambda,\phi\rangle_{\partial K_e}\Bigr)
  \\
  &+\sum_{K_e}\bigl\langle\boldsymbol v\cdot\boldsymbol \nu
    +\tau\,(p-\lambda),\ q\bigr\rangle_{\partial K_e}
     -\sum_{\mathfrak f^{\mathrm A}}c_e\,\langle\lambda,q\rangle_{\mathfrak f} .
\end{aligned}
\label{eq:calA_expanded}
\end{equation}
where $\mathfrak f^{\mathrm A}$ denotes the faces on $\Gamma_A$; Neumann faces carry no additional term.\\
Take  $y=Tx$, we get the energy identity
\begin{equation}
  \operatorname{Re}\,\mathcal{A}(x,Tx) = E(x) :=
  \sum_{K_e}\Bigl(s\,\rho\|\boldsymbol v\|_{K_e}^2
  + s\,\mu\,\|p\|_{K_e}^2
  + \tau\,\|p-\lambda\|_{\partial K_e}^2\Bigr)
   + \sum_{\mathfrak f^{\mathrm A}}c_e\,\|\lambda\|_{\mathfrak f}^2 ,
\label{eq:energy}
\end{equation}
where we used $\operatorname{Re}(-\sigma)=s>0$ and the fact that for any complex number $z\in\mathbb{C}$, we have $\operatorname{Re}(\bar{z} - z)= 0$. \\
Next, we bound the energy identity. It is straightforward that 
\begin{equation}
    \|\boldsymbol v\|_{L^2(\Omega)}^2\le(s\,\rho_{\min})^{-1}E(x)  \quad \text{and} \quad
\|p\|_{L^2(\Omega)}^2\le(s\mu_{\min})^{-1}E(x),
\end{equation}
For the face terms,  we have $\|\lambda\|_{\mathfrak f}^2
 \le 2\|p-\lambda\|_{\mathfrak f}^2+2\|p\|_{\mathfrak f}^2$
 and 
 \begin{equation}
  2\sum_{K}\|p-\lambda\|_{\partial K}^2
  \;\le\; \frac{2}{\tau_{\min}}\sum_{K}\tau\,\|p-\lambda\|_{\partial K}^2
  \;\le\; \frac{2}{\tau_{\min}}\,E(x) .
\label{eq:face_term1}
\end{equation}
and using
 the discrete trace inequality (see for example \cite{ern2021finite})
$\|q\|_{\partial K}^2\le C_{\mathrm{tr}}h_K^{-1}\|q\|_K^2$ on each cell with $C_{\mathrm{tr}}$ denoting the inequality constant,
followed by $\|p\|_\Omega^2\le E(x)/(s\,\mu_{\min})$, gives
\begin{equation}
  2\sum_{K}\|p\|_{\partial K}^2
  \;\le\; \frac{2\,C_{\mathrm{tr}}}{h_{\min}}\,\|p\|_\Omega^2
  \;\le\; \frac{2\,C_{\mathrm{tr}}}{h_{\min}\,s\,\mu_{\min}}\,E(x) .
\label{eq:face_term2}
\end{equation}
Finally we get 
\begin{equation}
  \sum_{\mathfrak f^{\mathrm I}}\|\lambda\|_{\mathfrak f}^2
  \;\le\; \Bigl(\frac{2}{\tau_{\min}}
        + \frac{2\,C_{\mathrm{tr}}}{h_{\min}\,s\,\mu_{\min}}\Bigr)E(x) .
\label{eq:face_final}
\end{equation}
While on ABC faces, with $c_{\min}:=\displaystyle \min_{\Gamma_A}(c\rho)^{-1}>0$, we have
\begin{equation}
\sum_{\mathfrak f^{\mathrm A}}\|\lambda\|_{\mathfrak f}^2
\le c_{\min}^{-1}E(x).
\end{equation}
Adding the five bounds gives
$\|x\|_{\mathcal{V}_h}^2\le\beta^{-1}E(x)\le \, \operatorname{Re}\,  \mathcal{A}(x,Tx) $ with $\beta$ defined as
\begin{equation}
  \beta^{-1} := \frac{1}{s\rho_{\min}}
  + \frac{1}{s\,\mu_{\min}}
  + \frac{2}{\tau_{\min}}
  + \frac{2\,C_{\mathrm{tr}}}{h_{\min}\,s\,\mu_{\min}}
  + \frac{1}{c_{\min}} ,
\label{eq:beta}
\end{equation}

\end{proof}
 
\begin{proposition}[Well posedness of the HDG system]
 Assume $s>0$, $\tau\ge\tau_{\min}>0$ and
$\mu\ge\mu_{\min}>0$. Then,
\begin{equation}
  \inf_{x\in \mathcal{V}_h\setminus\{0\}}\ \sup_{y\in \mathcal{V}_h\setminus\{0\}}
  \frac{|\mathcal{A}(x,y)|}{\|x\|_{\mathcal{V}_h}\,\|y\|_{\mathcal{V}_h}}\ \ge\ \beta ,
  \qquad
  \inf_{y\in \mathcal{V}_h\setminus\{0\}}\ \sup_{x\in \mathcal{V}_h\setminus\{0\}}
  \frac{|\mathcal{A}(x,y)|}{\|x\|_{\mathcal{V}_h}\,\|y\|_{\mathcal{V}_h}}\ \ge\ \beta .
\label{eq:infsup}
\end{equation}
In particular, the HDG system \eqref{eq:full_block} is well-posed. 
\end{proposition}

\begin{proof}
Fix $x\in \mathcal{V}_h$, since $Tx\in \mathcal{V}_h$ and $\|Tx\|_{\mathcal{V}_h}=\|x\|_{\mathcal{V}_h}$, the first inequality  in  \eqref{eq:infsup} is straightforward. 
For the second one, fix
$y\in \mathcal{V}_h\setminus\{0\}$ and we take $X=TY$.  We get 
$\mathcal{A}(Ty,y)=\mathcal{A}(x,Tx)\ge\beta\|x\|_{\mathcal{V}_h}^2=\beta\|y\|_{\mathcal{V}_h}^2$. Well-posedness follows from Banach-Nečas-Babuška theorem (see for instance \cite{fortin1991mixed}).
\end{proof}

\begin{proposition}[Local invertibility]
Under the same assumptions, $\mathbb{A}_e$ is coercive, hence invertible
for every cell, and the condensation   \eqref{eq:backsubstitution} is well defined.
 \end{proposition}

\begin{proof}
Let $u_e=(p,\boldsymbol v)\in W_h(K_e)\times\boldsymbol W_h(K_e)$  and apply the energy identity   \eqref{eq:energy} to
$x=(u_e,0)$, we get  $\mathcal{A}(x,Tx)=a_e(u_e,u_e)$ and  \eqref{eq:energy}
reduces to
\[
  \operatorname{Re}\,a_e(u_e,u_e)
  = s\,\rho\|\boldsymbol v\|_{K_e}^2
  + s\,\mu\,\|p\|_{K_e}^2
  + \tau\,\|p\|_{\partial K_e}^2
  \ \ge\ \alpha_a\|u_e\|_{L^2(K_e)}^2 ,
  \qquad \alpha_a := \min(s\,\rho,\,s\,\mu) > 0 ,
\]
so $\mathbb{A}_e$ is coercive and therefore invertible.
\end{proof}

The results we established are sufficient to prove well posedness of both the condensed and uncondnesed systems. Next , we establish coercivity of the sesquilinear   form $-\mathcal{K}$ that will be used for error estimation in Section \ref{sec:error_bounds}.

\begin{proposition}[Coercivity of $-\mathcal{K}$]
    Under the same assumptions, $-\mathcal{K}$ is coercive,
\begin{equation}
  \operatorname{Re}\bigl(-\mathcal{K}(\lambda,\lambda)\bigr)
  \ \ge\ \beta\,\|\lambda\|_{L^2(\Sigma)}^2
  \qquad\forall\,\lambda\in M_h.
\label{eq:Kcoerc}
\end{equation}
\end{proposition}

\begin{proof}
Let $\lambda\in M_h$ and $x_\lambda:=(u,\lambda)$ with
$U=-\mathbb{A}^{-1}\mathbb{C}\Lambda$, which yields $\mathcal{A}\bigl(x_\lambda,(u,0)\bigr)=0$, we obtain thus 
\[
  \mathcal{A}(x_\lambda,Tx_\lambda)
  = \mathcal{A}\bigl(x_\lambda,(u,0)\bigr) + \mathcal{A}\bigl(x_\lambda,(0,-\lambda)\bigr) 
  = \mathcal{A}\bigl(x_\lambda,(0,-\lambda)\bigr)  =  -\,\mathcal{K}(\lambda,\lambda),
\]
which implies that 
\begin{equation}
\label{eq:coercivity_K}
    \operatorname{Re}(-\mathcal{K}(\lambda,\lambda))
=\operatorname{Re}\mathcal{A}(x_\lambda,Tx_\lambda)
\ge\beta\|x_\lambda\|_{\mathcal{V}_h}^2\ge\beta\|\lambda\|_{L^2(\Sigma)}^2.
\end{equation}
\end{proof}

\begin{remark}
The coercivity constant $\beta_k$ of $-\mathcal{K}$ satisfies  $\beta_k>\beta$ in practice, indeed
the constant $\beta$ in \eqref{eq:coercivity_K} must hold 
for all variables $x_\lambda=(u,\lambda)$, whereas the coercivity of
$-\mathcal{K}$ only needs to control $\lambda$. Indeed, from \eqref{eq:coercivity_K}, we have
\begin{equation}
     \operatorname{Re}(-\mathcal{K}(\lambda,\lambda))\ge\beta\|x_\lambda\|_{\mathcal{V}_h}^2  = \beta\|\lambda\|_{L^2(\Sigma)}^2 + \beta\|u_\lambda\|_{L^2(\Omega)}^2.
\end{equation}

\end{remark}

\begin{remark}[Choice of the stabilization parameter]
The analysis we provided only requires $\tau\ge\tau_{\min}>0$. Following
\cite{faucher2020adjoint}, we take a media-dependent stabilization parameter
\begin{equation}
   \tau := \rho^{-1},
\label{eq:tau_def}
\end{equation}
which is the diffusion coefficient of the second-order form
$-\nabla\cdot(\rho^{-1}\nabla p)+\sigma^{2}\mu\,p=-\sigma f$, the scaling was proposed for linear convection--diffusion equations in \cite{nguyen2009implicit}.
\end{remark}

\subsection{Analyticity of the solution map}
\label{sec:analyticity}
In this section, we establish the analyticity of the HDG solution map 
\begin{equation}
    \sigma \in \mathcal{D}_\sigma \mapsto X(\sigma) := (U(\sigma), \Lambda(\sigma)),
    \label{eq:HDG_sol_map}
\end{equation}
where $X(\sigma)$ denotes the coefficient vector of the solution of  \eqref{eq:full_block} at $\sigma$ and $\mathcal{D}_\sigma := \{\, i\omega - s \;:\; \omega\in\mathcal{D}_\omega \}$, where $s>0$ is given and $\mathcal{D}_\omega \subset \mathbb{R}^*_+$ a frequency range.\\
The analyticity allows the enrichment of the reduced basis with frequency derivatives in Section \ref{sec:basis_construction}. We recall that the continuous solution map 
\begin{equation}
    \sigma \mapsto (p(\sigma), \boldsymbol{v}(\sigma)),
\end{equation}
is analytic since the PDE is governed by a linear Fredholm operator that is invertible, thus the inverse is also analytic. 
For the HDG solution map \eqref{eq:HDG_sol_map}, we know that matrix $\mathbb{A}$ is a matrix polynomial of
degree $1$ with respect to $\sigma$ and is invertible (local solvability). Cramer's rule gives
\[
  \mathbb{A}(\sigma)^{-1}
  =
  \frac{\operatorname{adj}(\mathbb{A}(\sigma))}{\det(\mathbb{A}(\sigma))}.
\]
Hence, $\sigma\mapsto \mathbb{A}(\sigma)^{-1}$ is analytic on $\mathcal{D}_\sigma $. 
With $\mathbb{A}(\sigma)^{-1}$ analytic, so is 
\begin{equation}
  \mathbb{K}(\sigma) \;:=\; \mathbb{L} - \mathbb{B}\,\mathbb{A}(\sigma)^{-1}\,\mathbb{C}. 
  \label{eq:schur}
\end{equation}
Similarly, since $\mathbb{K}(\sigma)$ is invertible on
$\mathcal{D}_\sigma$ (global solvability), we get that  its inverse
$\mathbb{K}(\sigma)^{-1}$ is analytic on $\mathcal{D}_\sigma$.  The skeleton solution $\Lambda(\sigma)$  is thus analytical and by substitution so is $U(\sigma)$.

\section{Reduced basis approximation}
\label{sec:ROM}
Reduced basis methods seek a low-rank approximation of the high-fidelity HDG solutions. Although high-fidelity solutions live in the high-dimensional space $\mathcal{V}_h$, the set of all
solutions obtained as the frequency $\sigma$ varies over a given domain $\mathcal{D}_\sigma$, 
\[
\mathcal{S}_h := \{(u_h(\sigma),\lambda_h(\sigma)) : \sigma\in \mathcal{D}_\sigma\},
\]
can be accurately approximated by a linear space of much smaller
dimension. Indeed, solutions corresponding to nearby frequencies typically exhibit similar features, so that a relatively small number of high-fidelity solutions, referred to as \emph{snapshots}, can contain most of the information required to approximate the solution at other frequencies in $\mathcal{D}_\sigma$.

A common tool for quantifying the approximability of $\mathcal{S}_h$ by finite-dimensional linear spaces is the Kolmogorov $n$-width,
\begin{equation}
\label{eq:nwidth}
    d_n(\mathcal{S}_h, \mathcal{V}_h) := \inf_{\substack{\mathcal{Z} \subset \mathcal{V}_h \\ \dim(\mathcal{Z}) = n}} \; \sup_{s \in \mathcal{S}_h} \; \inf_{v \in \mathcal{Z}} \, \| s - v \|_{\mathcal{V}_h}.
\end{equation}
Introduced by Kolmogorov~\cite{kolmogoroff1936beste}, this quantity measures the best possible worst-case approximation error obtained by $n$-dimensional subspaces of $\mathcal{V}_h$. In the context of time-harmonic wave propagation, exponential decay of the Kolmogorov $n$-width has been established in \cite{henriquez2026fast,henriquez2026model}, in contrast to the algebraic decay observed for time-domain wave propagation~\cite{greif2019decay}. Nevertheless, even in the time-harmonic setting, the reduced dimension may become relatively large before the asymptotic exponential decay regime is reached, particularly for highly oscillatory solutions. In such cases, one may resort to truncation strategies to control the reduced dimension, or consider nonlinear model reduction approaches. One such approach consists in decomposing the parametric domain into smaller subsets and constructing a local reduced model for each subset, thereby allowing to control the rank of the reduced spaces ~\cite{tifouti2025registration,eftang2010hp,Haasdonk}.

The optimization problem \cref{eq:nwidth} defining the Kolmogorov $n$-width is inherently nonlinear and non-convex. Although the existence of optimal subspaces has been established in reflexive Banach spaces~\cite{garkavi1962optimal}, constructing such optimal spaces is generally impractical. Model reduction methods therefore aim to construct computationally tractable reduced spaces that provide a good approximation of the optimal subspaces.

Reduced basis methods for parametric problems  exploit an \emph{offline}/\emph{online} computational framework:  in an
\emph{offline} stage, performed once, the high-fidelity problem is solved at
the selected frequencies and the snapshots are used to build a low-dimensional
basis, and the FOM
is projected onto this basis  resulting in a small system. In an \emph{online} stage, for any new frequency $\sigma\in \mathcal{D}_\sigma$ we compute the \emph{reduced solution} $(u_r,\lambda_r)$. The offline cost is thus compensated over many online evaluations, for instance
a dense frequency sweep.

More precisely, we seek reduced approximations  
\[u_r\in\mathcal{Z}_n^U:=\operatorname{span}(\boldsymbol{Z}_U) \quad \text{and} \quad  \lambda_r\in\mathcal{Z}_n^\Lambda:=\operatorname{span}(\boldsymbol{Z}_\Lambda)\] of ranks
$n_U\ll N_U$, $n_\Lambda\ll N_\Lambda$.  Note that the reduced order
basis $\bs{Z}_U=(\zeta^U_1,\ldots,\zeta^U_{n_U}) \subset W_h\times\boldsymbol W_h$ and
$\bs{Z}_\Lambda=(\zeta^\Lambda_1,\ldots,\zeta^\Lambda_{n_\Lambda}) \subset M_h$, and we define the coupled reduced space $\mathcal{Z}_n:=\mathcal{Z}_n^U\times\mathcal{Z}_n^\Lambda$.\\\ 

We further denote by
$\mathbb{Z}_U\in\mathbb{C}^{N_U\times n_U}$ and
$\mathbb{Z}_\Lambda\in\mathbb{C}^{N_\Lambda\times n_\Lambda}$ the associated
coefficient vectors of the basis functions; note that the bases are orthonormal in $L^2(\Omega)$ and $L^2(\Sigma)$ i.e 
$\mathbb{Z}_U^{*}\mathbb{M}_U\mathbb{Z}_U=I$ and
$\mathbb{Z}_\Lambda^{*}\mathbb{M}_\Lambda\mathbb{Z}_\Lambda=I$.\\
The reduced solutions are defined using the linear ansatz
\begin{equation}
   u_r := \bs{Z}_U\,\widehat{U},
    \qquad
     \lambda_r := \bs{Z}_\Lambda\,\widehat{\Lambda},
    \label{eq:reduced_ansatz}
\end{equation}
where $\widehat{U}\in \mathbb{C}^{n_U}$ and $\widehat{\Lambda}\in\mathbb{C}^{n_\Lambda}$ are the corresponding generalized coefficients.

We recall that since the ABC boundary condition  in \eqref{eq:acoustic_strong}
is expressed directly in the physical fields and is therefore
frequency-independent, the uncondensed operator  \eqref{eq:full_block}
is affine in $\sigma$: we define
\begin{equation}
   \mathbb{G}(\sigma) := \mathbb{G}_0 + \sigma\,\mathbb{G}_1,
   \qquad
   \mathbb{G}_0 :=
   \begin{pmatrix} \mathbb{A}_0 & \mathbb{C}\\ \mathbb{B} & \mathbb{L}\end{pmatrix},
   \qquad
   \mathbb{G}_1 :=
   \begin{pmatrix} \mathbb{A}_1 & 0\\ 0 & 0\end{pmatrix},
   \label{eq:G_affine}
\end{equation}
where $\mathbb{A}=\mathbb{A}_0+\sigma\mathbb{A}_1$ splits  \eqref{eq:Ae_def}
and the right-hand side $\mathbb{R}_S:=(\mathbb{S}, 0)^T$ is assumed frequency-independent. \\
Building the reduced model follows two steps, first we build the reduced order basis (ROB) in Section \ref{sec:basis_construction}, second, we project the full model \eqref{eq:full_block} onto the reduced basis in Section \ref{sec:galerkin}. 

\subsection{Construction of the reduced basis}
\label{sec:basis_construction}
The bases $(\bs{Z}_U,\bs{Z}_\Lambda)$ are built using a greedy loop over the frequency range. 
Each iteration of the greedy algorithm selects the least well approximated frequency by means of a residual indicator, computes the solution at that frequency together with a few of its frequency derivatives and augments them to the snapshot set, which is then compressed by a weighted POD. The three
ingredients --- moment enrichment, the residual indicator, and the
weighted POD compression --- are detailed below and combined in the
overall procedure of \cref{alg:overall} in Section  \ref{sec:overall_algo}.

\subsubsection{Enriching the reduced basis with frequency derivatives} 

At each selected training  point $\sigma_{\ell}$,  we enrich the snapshot set
with the solution and its frequency derivatives. Indeed, rather than
representing the parametric solution by a single snapshot per selected
frequency, we exploit the smoothness of the map $\sigma\mapsto X(\sigma)=(U(\sigma), \Lambda(\sigma))$ (see Section \ref{sec:analyticity}) and
append a few of its frequency derivatives at $\sigma_{\ell}$, so that local slope information is also added to the training set.

To get the analytical derivatives of the solution we expand the solution
\begin{equation}
    X(\sigma) = \sum_{j\ge0}(\sigma-\sigma_{\ell})^j X^{(j)}
    \label{eq:X_taylor}
\end{equation}
and the PDE operators
\begin{equation}
    \mathbb{G}(\sigma) = \sum_{j\ge0}(\sigma-\sigma_{\ell})^j  \mathbb{G}^{(j)}
    \label{eq:K_taylor}
\end{equation}
Then we identify powers of $(\sigma-\sigma_{\ell})$ in the expanded uncondensed
system  \eqref{eq:full_block}. We get the multi-RHS system
\begin{equation}
    \mathbb{G}^{(0)}X^{(0)} = \mathbb{R}_S,
    \qquad
    \mathbb{G}^{(0)}X^{(j)} = -\,\mathbb{G}^{(1)}X^{(j-1)},
    \quad j\ge1,
    \label{eq:taylor_hierarchy}
\end{equation}
with $\mathbb{G}^{(0)}=\mathbb{G}(\sigma_{\ell})$ and
$\mathbb{G}^{(j)}=\frac{1}{j!}\,\partial^{j}_\sigma\mathbb{G}(\sigma_{\ell})$. Since  $\mathbb{G}$ is affine in $\sigma$, all derivatives beyond
the first vanish, $\partial_\sigma^{j}\mathbb{G}=0$ for $j\ge2$. In practice, we follow the Schur elimination presented in Section \ref{sec:condensed} to get the derivatives --  since only the local blocks
$\mathbb{A}_e$ and the condensed operator $\mathbb{K}$ are factorized in practice. \\
Each derivative is recovered by
static condensation of its right-hand side,
\begin{equation}
   \mathbb{K}\,\Lambda^{(j)}
   = \sum_e \mathcal{R}_e^{T}\,\mathbb{B}_e\,\mathbb{A}_e^{-1}\,
             \mathbb{A}_{1,e}\,U_e^{(j-1)} ,
   \qquad
   U_e^{(j)} = \mathbb{A}_e^{-1}\Bigl(-\mathbb{A}_{1,e}U_e^{(j-1)}
              - \mathbb{C}_e\,\mathcal{R}_e\Lambda^{(j)}\Bigr),
   \label{eq:moments_condensed}
\end{equation}
All the derivatives
$\{X^{(j)}\}_{j\geq0}$ therefore reuse the same factorizations and are
obtained at the cost of one additional back-substitution each, making
the enrichment step cheap relative to the high-fidelity solve. Further details on the HPC implementation and computational aspects are provided in Section~\ref{sec:overall_algo}.

\subsubsection{Greedy selection of training points}
\label{sec:greedy_selection}

In the weak greedy loop, we use the well-known \emph{a posteriori} residual-based indicator to select the expansion points. In Section \ref{sec:wellposedness}, we proved that the discrete system \eqref{eq:full_block} is inf-sup stable; in this case the residual-based indicator is known to reproduce the behaviour of the actual approximation error \cite{binev2011convergence}, without requiring the computation of the high-fidelity solutions.
Let $\mathcal{D}^h_\sigma\subset\mathcal{D}_\sigma$ be a finite training
set. Starting from an initial ROB $(\mathbb{Z}_U, \mathbb{Z}_\Lambda)$, the next step of the greedy algorithm enriches the ROB with the training point that maximizes
\begin{equation}
    \Delta(\sigma) := \frac{\|r(\sigma)\|_{\mathcal{V}_h'}}{\|\mathbb{R}_S\|_{\mathcal{V}_h'}},
    \qquad
    r(\sigma) := \mathbb{G}(\sigma)\,(\mathbb{Z}_U\widehat{U}, \mathbb{Z}_\Lambda\widehat{\Lambda})^T - \mathbb{R}_S, 
    \label{eq:greedy_indicator}
\end{equation}
where  $\|.\|_{\mathcal{V}_h'}$ is the dual norm of $\|\cdot\|_{\mathcal{V}_h}$, $\|r\|_{\mathcal{V}_h'}^2 = r^{*}\,\mathbb{M}^{-1} r$.\\
The indicator quantifies,
at a cost independent of the full-order solve, how well the current reduced
solution satisfies the full-order equations at a candidate frequency, and
thus serves as a computable surrogate for the true error in the greedy
selection.
\paragraph{Assembling the residual-based indicator. }
We note that the indicator \eqref{eq:greedy_indicator} is never globally assembled: 
 since $\mathbb{G}$ is affine, the residual $\|r(\sigma)\|_{\mathcal{V}_h'}^2$ is a quadratic polynomial in
$\sigma$ and $\bar\sigma$, so it is precomputed once per basis update through element-wise matrix multiplications of
Gram matrices in the reduced coordinates, and evaluated online at cost
$O((n_U+n_\Lambda)^2)$ per candidate frequency, independently of the high-fidelity order $N$. This assembly of the residual indicator is attractive for the greedy loop since the basis is updated at each iteration and the indicator is evaluated over hundreds of frequency points to select the least well-approximated one, for which a raw evaluation would be prohibitive. Instead, the assembly follows an assembly--evaluation splitting.\\
\medskip\noindent\emph{1. Assembly.}\\
With $(u_r,\lambda_r)=(\,\bs{Z}_U\hat{U},\ \bs{Z}_\Lambda\hat{\Lambda}\,)$, the two residual blocks read
\begin{equation}
\begin{aligned}
    r_U(\sigma) &= (\mathbb{A}_0+\sigma\mathbb{A}_1)\mathbb{Z}_U\hat{U}
                + \mathbb{C}\,\mathbb{Z}_\Lambda\hat{\Lambda} - \mathbb{S},\\
    r_\Lambda(\sigma) &= \mathbb{B}\,\mathbb{Z}_U\hat{U}
                + \mathbb{L}\,\mathbb{Z}_\Lambda\hat{\Lambda} .
\end{aligned}
    \label{eq:residual_blocks}
\end{equation}
The projected operators $\mathbb{A}_0\mathbb{Z}_U$, $\mathbb{A}_1\mathbb{Z}_U$, $\mathbb{C}\mathbb{Z}_\Lambda$, $\mathbb{B}\mathbb{Z}_U$, $\mathbb{L}\mathbb{Z}_\Lambda$ and the source $\mathbb{S}$ are computed element-wise, then multiplied by the sparse norm matrix $\mathbb{M}^{-1}$ into reduced Gram matrices and vectors, once per basis update. For instance, the volume Gram matrices read
\begin{equation}
    H_{00} := (\mathbb{A}_0\mathbb{Z}_U)^{*} \mathbb{M}_U^{-1} (\mathbb{A}_0\mathbb{Z}_U),
    \quad
    H_{01} := (\mathbb{A}_0\mathbb{Z}_U)^{*} \mathbb{M}_U^{-1} (\mathbb{A}_1\mathbb{Z}_U),
    \quad
    H_{11} := (\mathbb{A}_1\mathbb{Z}_U)^{*} \mathbb{M}_U^{-1} (\mathbb{A}_1\mathbb{Z}_U).
    \label{eq:gram_volume}
\end{equation}
The same logic holds for the analogous skeleton and cross-block Gram matrices, as well as for the source term. We note that the Gram matrices are of size at most $(n_U+n_\Lambda)\times(n_U+n_\Lambda)$, which makes the evaluation of the residual over many frequencies affordable in the greedy selection. 

\medskip\noindent\emph{2. Evaluation.}\\
At each candidate frequency, $\|r(\sigma)\|_{\mathcal{V}_h'}^2$ is recovered blockwise using the identity
\begin{equation}
     \|r_U(\sigma)\|_{\mathcal{V}_h'}^2 = r_U(\sigma)^*    \mathbb{M}_U^{-1}    r_U(\sigma),
\end{equation}
so that the volume block yields the scalar quadratic form
\begin{equation}
    \|r_U(\sigma)\|_{\mathcal{V}_h'}^2
    = \hat{U}^{*} H_{00}\hat{U}
    + 2\,\mathrm{Re}\!\bigl(\sigma\,\hat{U}^{*} H_{01}\hat{U}\bigr)
    + |\sigma|^2\,\hat{U}^{*} H_{11}\hat{U}
    + (\text{cross and source terms}),
    \label{eq:res_quadratic}
\end{equation}
where $H_{00},H_{01},H_{11}$ are constant per basis update. The same holds for $\|r_\Lambda(\sigma)\|_{\mathcal{V}_h'}^2$, which gives 
\begin{equation}
    \Delta(\sigma) := \frac{\sqrt{ \|r_U(\sigma)\|_{\mathcal{V}_h'}^2 +  \|r_\Lambda(\sigma)\|_{\mathcal{V}_h'}^2}}{\|\mathbb{R}_S\|_{\mathcal{V}_h'}}.
\end{equation}
Each evaluation at $\sigma$ scales as $O((n_U+n_\Lambda)^2)$, which allows the training set $\mathcal{D}^h_\sigma$ to be swept at negligible cost relative to a single high-fidelity solve. We provide further HPC implementation and computational aspects in Section~\ref{sec:overall_algo}. 

\medskip\noindent
The greedy loop terminates when the true relative error at the frequency
selected by the residual indicator, that is, the relative $L^2$ error between
the high-fidelity solution and the reduced solution at that frequency, falls
below the prescribed tolerance. The indicator itself is not compared with this
threshold, since it differs from the relative error by a scaling factor.

\subsubsection{POD compression}
The snapshots and moments accumulated over the greedy iterations are
compressed by POD and yield orthonormal bases. The volume and skeleton snapshots $\{(u_h,\lambda_h)(\sigma)\}_\sigma$ live on different
meshes and are compressed separately, using the $L^2(\Omega)$ and
$L^2(\Sigma)$ inner products rather than the Euclidean one. This ensures that the
retained modes are optimal in their respective variational norms :
\begin{equation}
    \bs{Z}_U = \mathrm{POD}_{L^2(\Omega)}\bigl\{u_h(\sigma_j)\bigr\}_j,
    \qquad
    \bs{Z}_\Lambda = \mathrm{POD}_{L^2(\Sigma)}\bigl\{\lambda_h(\sigma_j)\bigr\}_j.
    \label{eq:ROB_extraction}
\end{equation}
Compressing the two fields independently, rather than stacking them into
a single snapshot, keeps the two reduced dimensions $n_U$ and $n_\Lambda$
decoupled and preserves the block structure of  \eqref{eq:full_block} at
the reduced level, as exploited in Section  \ref{sec:galerkin}.

\subsection{Galerkin projection}
\label{sec:galerkin}
In this section, we construct the reduced order model by Galerkin projection of the non-condensed system \eqref{eq:full_block}. As noted in \cref{sec:condensed}, projecting the condensed system \eqref{eq:condensed} directly destroys the affine parameter dependence that is explicit in the non-condensed HDG formulation which allows an efficient offline/online decomposition. Moreover, even if the reduced condensed system were made online-efficient through hyper-reduction, the volume unknowns would still have to be recovered from the hybrid variable by local back-substitution. Although this reconstruction is local, it must be carried out on every mesh element, and therefore adds an online cost that scales with the size of the high-fidelity discretization rather than with the dimension of the reduced space. 

Once the ROBs $\bs{Z}_U$ and $\bs{Z}_\Lambda$ are available, we proceed to building the ROM.
We exploit the structure of the system presented in Section \ref{sec:uncondensed} and projecting   \eqref{eq:full_block} onto $\boldsymbol{Z}_U \times\boldsymbol{Z}_\Lambda$ gives the reduced model. For each query point $\omega \in \mathcal{D}_\omega$, we set $\sigma= i\omega -s$, and solve: 
\begin{equation}
    \begin{pmatrix}
    \widehat{\mathbb{A}}(\sigma) & \widehat{\mathbb{C}} \\[4pt]
    \widehat{\mathbb{B}} & \widehat{\mathbb{L}}
    \end{pmatrix}
    \begin{pmatrix} \widehat{U} \\ \widehat{\Lambda} \end{pmatrix}
    =
    \begin{pmatrix} \widehat{\mathbb{S}} \\ 0 \end{pmatrix},
    \label{eq:reduced_block}
\end{equation}
with
\begin{equation}
    \widehat{\mathbb{A}}(\sigma) := \mathbb{Z}_U^*\,\mathbb{A}(\sigma)\,\mathbb{Z}_U,
    \quad
    \widehat{\mathbb{C}} := \mathbb{Z}_U^*\,\mathbb{C}\,\mathbb{Z}_\Lambda,
    \quad
    \widehat{\mathbb{B}} := \mathbb{Z}_\Lambda^*\,\mathbb{B}\,\mathbb{Z}_U,
    \quad
    \widehat{\mathbb{L}} := \mathbb{Z}_\Lambda^*\,\mathbb{L}\,\mathbb{Z}_\Lambda,
    \quad
    \widehat{\mathbb{S}} := \mathbb{Z}_U^*\,\mathbb{S},
\label{eq:reduced_operators}
\end{equation}
of size $(n_U+n_\Lambda)\ll(N_U+N_\Lambda)$, and preserving the block
structure of  \eqref{eq:full_block}.\\
As already stated in Remark \ref{rem:local_proj}, we never proceed to assembling the global matrices, instead and thanks to the HDG element-wise structure, we perform a local element-wise projection of operators, thus enhancing memory and CPU cost performances. Further implementation details are provided in Section~\ref{sec:overall_algo}.

This projection is performed once, in an offline stage whose cost scales with the full-order dimension $N$, and its outcome is reused for every online query. Indeed, thanks to the affine frequency dependency, using the splitting  \eqref{eq:G_affine} the reduced operators are built offline and stored. Once the query frequency is available online, we only assemble the cheap reduced model online. 

\paragraph{Online evaluation.}
For each query frequency $\sigma\in\mathcal{D}_\sigma$, the $\sigma$-dependent block is reassembled from the two stored matrices as $\widehat{\mathbb{A}}(\sigma) = \widehat{\mathbb{A}}_0 + \sigma\,\widehat{\mathbb{A}}_1$, at cost $O(n_U^2)$ --- a scalar combination of two small dense matrices. Eliminating the volume unknown by a reduced Schur complement then yields the reduced model
\begin{equation}
    \widehat{\mathbb{K}}(\sigma)\,\widehat{\Lambda}
    = \widehat{\mathbb{F}}(\sigma),
        \label{eq:condensed_reduced}
\end{equation}
where
\begin{equation}
    \widehat{\mathbb{K}}(\sigma)
    := \widehat{\mathbb{L}}
     - \widehat{\mathbb{B}}\,\widehat{\mathbb{A}}(\sigma)^{-1}\widehat{\mathbb{C}}, 
    \qquad \text{and}   \quad     
    \widehat{\mathbb{F}}(\sigma)
    := -\widehat{\mathbb{B}}\,\widehat{\mathbb{A}}(\sigma)^{-1}\widehat{\mathbb{S}}. 
\end{equation}
The volume unknown is then recovered 
\begin{equation}
    \widehat{U} = \widehat{\mathbb{A}}(\sigma)^{-1}
(\widehat{\mathbb{S}} - \widehat{\mathbb{C}}\widehat{\Lambda}).
\label{eq:reduced_backsub}
\end{equation}
Note that the online cost of the reduced model is independent of $N$ and scales as $O(n_U^3)$ for the inversion of $\widehat{\mathbb{A}}(\sigma)$ and $O(n_\Lambda^3)$ for the solution of the reduced Schur complement system \eqref{eq:condensed_reduced}.

\begin{remark}
\label{rem:wellposdeness_ROM}
The well-posedness of the reduced model  \eqref{eq:reduced_block} is inherited from
Section  \ref{sec:wellposedness}. Since the two bases are constructed separately, the reduced
space $\mathcal{Z}_n=\mathcal{Z}_n^U\times\mathcal{Z}_n^\Lambda$ is $T$-invariant: for
$x_r=(u_r,\lambda_r)\in\mathcal{Z}_n$ we have $Tx_r=(u_r,-\lambda_r)\in\mathcal{Z}_n$, so
$Tx_r$ is an admissible test function and the energy identity  \eqref{eq:energy} still applies 
on $\mathcal{Z}_n$. Hence $\operatorname{Re}\mathcal{A}(x_r,Tx_r)\ge\beta\|x_r\|^2_{\mathcal{V}_h}$
with the same constant $\beta$ as the full-order system. In contrast, using  a single basis stacking the volume and skeleton components, $Tx_r$ would
in general leave $\mathcal{Z}_n$ and the argument breaks. 
Moreover, the coercivity of the sesquilinear   form $a$ defined in \eqref{eq:form_ae} justifies the invertibility of
the reduced matrix in  \eqref{eq:reduced_block} which in turn justifies \cref{eq:reduced_backsub}. Finally, since $\beta$ is the inf-sup constant of the full-order system, it is independent of the reduced dimension $n_U$ and $n_\Lambda$ : meaning that the reduced model remains uniformly stable as the reduced bases are enriched. 
\end{remark}

\subsection{Error analysis}
\label{sec:error_bounds}
We provide in this section an a priori error analysis of the proposed reduced model \cref{eq:condensed_reduced,eq:reduced_backsub}. First, we derive a general error estimate that couples both volume and skeleton variables using the classical C\'ea-Babuška estimate. Second, we provide error estimates that decouple the two variables and provides a systematic way of enriching the basis with a lifting method.

We introduce the errors between high-fidelity solutions and reduced solutions for both unkowns:
\begin{equation}
    e_U:=u_h-u_r
\end{equation} and  
\begin{equation}
    e_\Lambda:=\lambda_h-\lambda_r
\end{equation}

\subsubsection{Coupled error bound}
The following result shows that the joint reduced solution $x_r=(u_r,\lambda_r)\in \mathcal{Z}_n$ is quasi-optimal: the error is bounded by the best approximation of the high-fidelity
solution in the reduced space $\mathcal{Z}_n$.

\begin{proposition}[C\'ea--Babuška lemma]\label{prop:cea}
Let $\gamma$ denote the continuity constant of the sesquilinear form
$\mathcal{A}$ and $\beta$ its inf-sup constant defined in \eqref{eq:beta}. Then
\begin{equation}\label{eq:cea}
  \bigl\|(e_U,e_\Lambda)\bigr\|_{\mathcal{V}_h}
  \le \Bigl(1+\frac{\gamma}{\beta}\Bigr)
  \inf_{(v_r,q_r)\in\mathcal{Z}_n}
  \bigl\|(u_h-v_r,\,\lambda_h-q_r)\bigr\|_{\mathcal{V}_h}.
\end{equation}
\end{proposition}

\begin{proof}
For $(v_r,q_r)\in\mathcal Z_n$ set $d_U:=u_h-v_r$,
$d_\Lambda:=\lambda_h-q_r$, $\zeta_U:=u_r-v_r$, $\zeta_\Lambda:=\lambda_r-q_r$,
so that we get  $(\zeta_U,\zeta_\Lambda)\in\mathcal Z_n$ and
$(e_U,e_\Lambda)=(d_U,d_\Lambda)-(\zeta_U,\zeta_\Lambda)$.\\ For any
$(w_r,\varphi_r)\in\mathcal Z_n$, by linearity and since $(u_h,\lambda_h)$
solves  \eqref{eq:full_block} on $\mathcal Z_n$ while $(u_r,\lambda_r)$
solves  \eqref{eq:reduced_block},
\[
  \mathcal A\bigl((\zeta_U,\zeta_\Lambda),(w_r,\varphi_r)\bigr)
  =\mathcal A\bigl((u_h,\lambda_h),(w_r,\varphi_r)\bigr)
   -\mathcal A\bigl((v_r,q_r),(w_r,\varphi_r)\bigr)
  =\mathcal A\bigl((d_U,d_\Lambda),(w_r,\varphi_r)\bigr) .
\]
Inf-sup stability and continuity of $\mathcal A$ then give
\[
  \beta\|(\zeta_U,\zeta_\Lambda)\|_{\mathcal{V}_h}
  \le\sup_{(w_r,\varphi_r)\in\mathcal Z_n\setminus\{0\}}
  \frac{\mathcal A((\zeta_U,\zeta_\Lambda),(w_r,\varphi_r))}{\|(w_r,\varphi_r)\|_{\mathcal{V}_h}}
  \le\gamma\|(d_U,d_\Lambda)\|_{\mathcal{V}_h} ,
\]
so $\|(\zeta_U,\zeta_\Lambda)\|_{\mathcal{V}_h}\le(\gamma/\beta)\|(d_U,d_\Lambda)\|_{\mathcal{V}_h}$.
We use the triangle inequality and the infimum over $(v_r,q_r)$ to conclude.
\end{proof}

\subsubsection{Decoupled error bounds}
Estimate  \eqref{eq:cea} bounds the joint error $\|(e_U,e_\Lambda)\|_{\mathcal{V}_h}$, but does not indicate  how well
$u_r$ and $\lambda_r$ separately approximate $u_h$ and $\lambda_h$. In this section, we derive decoupled error estimates for each unknown. We note that classical
saddle-point analysis for MOR (see for instance \cite{gerner2012certified}) does not apply
here, since $\mathbb{C}\neq\mathbb{B}^{T}$ and $\mathbb{L}\neq0$. Instead, the decoupled
estimates below follow from the coercivity of $\mathbb{A}$ and $-\mathbb{K}$
(see Section  \ref{sec:wellposedness}) and the Schur complement structure described in
Section  \ref{sec:condensed}. 
 
\paragraph{1. Preliminaries.}
We recall that the reduced model \eqref{eq:condensed_reduced} for $\lambda$ is obtained by projecting the
non-condensed system of Section~\ref{sec:uncondensed}. In contrast, the high-fidelity trace $\lambda_h$
is obtained by solving the condensed system of Section~\ref{sec:condensed}. As
already discussed there, we chose the non-condensed formulation because it
preserves the frequency-affine structure of the reduced problem and, unlike
the condensed formulation, allows us to build reduced models for both $u$ and
$\lambda$. The price to pay is that $\lambda_h$ and $\lambda_r$ no longer solve
the same reduced problem, since projection and Schur condensation do not
commute in general. We now derive error estimates that quantify this gap.

Projecting the condensed system 
 onto $\mathcal{Z}^\Lambda_n$ would give the following reduced operator and reduced
right-hand side
\[
   \mathbb{Z}_\Lambda^{*}\mathbb{K}\mathbb{Z}_\Lambda
   =\widehat{\mathbb{L}}
    -\mathbb{Z}_\Lambda^{*}\mathbb{B}\,\mathbb{A}^{-1}\,\mathbb{C}\mathbb{Z}_\Lambda ,
   \qquad
   \mathbb{Z}_\Lambda^{*}\mathbb{F}
   =-\mathbb{Z}_\Lambda^{*}\mathbb{B}\,\mathbb{A}^{-1}\,\mathbb{S}.
\]
The ROM \cref{eq:condensed_reduced} does the opposite: it projects first and condenses afterwards,
which by  \eqref{eq:reduced_operators} gives
\[
   \widehat{\mathbb{K}}
   =\widehat{\mathbb{L}}
    -\mathbb{Z}_\Lambda^{*}\mathbb{B}\,
     \bigl[\mathbb{Z}_U\widehat{\mathbb{A}}^{-1}\mathbb{Z}_U^{*}\bigr]\,
     \mathbb{C}\mathbb{Z}_\Lambda ,
   \qquad
   \widehat{\mathbb{F}}
   =-\mathbb{Z}_\Lambda^{*}\mathbb{B}\,
     \bigl[\mathbb{Z}_U\widehat{\mathbb{A}}^{-1}\mathbb{Z}_U^{*}\bigr]\,\mathbb{S}.
\]
Note that the difference between the two reduced models come from the way the volume unknown is eliminated. Indeed, by subtracting we get 
\begin{equation}
   \delta\mathbb{K}
   :=\mathbb{Z}_\Lambda^{*}\mathbb{K}\mathbb{Z}_\Lambda-\widehat{\mathbb{K}}
   =-\,\mathbb{Z}_\Lambda^{*}\mathbb{B}
      \bigl(\mathbb{A}^{-1}-\mathcal{P}\bigr)\mathbb{C}\mathbb{Z}_\Lambda ,
   \qquad
   \delta\mathbb{F}
   :=\mathbb{Z}_\Lambda^{*}\mathbb{F}-\widehat{\mathbb{F}}
   =-\,\mathbb{Z}_\Lambda^{*}\mathbb{B}
      \bigl(\mathbb{A}^{-1}-\mathcal{P}\bigr)\mathbb{S}.
   \label{eq:commutation}
\end{equation}
where 
\begin{equation}
   \mathcal{P}:=\mathbb{Z}_U\widehat{\mathbb{A}}^{-1}\mathbb{Z}_U^{*}.
   \label{eq:local_solve_op}
\end{equation}
 \medskip\noindent Note that $\mathcal{P}$ is a Galerkin projection. Indeed, let $f$ be a volume right-handside.  We have
 \begin{equation}  
 a(\mathcal{P}f, v_r)= \hat{v}^{*}\mathbb{Z}_U^{*}\mathbb{A}\mathbb{Z}_U\widehat{\mathbb{A}}^{-1}\mathbb{Z}_U^{*}f= \hat{v}^{*}\widehat{\mathbb{A}}\widehat{\mathbb{A}}^{-1}\mathbb{Z}_U^{*}f = v_r^{*} f = a(\mathbb{A}^{-1}f, v_r) .
 \end{equation}
 This proves that  $\mathcal{P}f$ is a ROM-Galerkin approximation of $\mathbb{A}^{-1}f$ in $\mathcal{Z}_n^U$, and we have by C\'ea's lemma 
 \begin{equation}   \bigl\|\mathbb{A}^{-1}f-\mathcal{P}f\bigr\|_{L^2(\Omega)}   \;\le\;\frac{\gamma_a}{\alpha_a}\,   \inf_{v_r\in\mathcal{Z}^U_n}     \bigl\|\mathbb{A}^{-1}f-v_r\bigr\|_{L^2(\Omega)} .   \label{eq:cea_local}
 \end{equation}

\paragraph{2. Trace estimate.} 
We first bound the trace error $e_\Lambda$. Proposition~\ref{prop:trace} shows that the bound
combines two terms: the best-approximation error, and how far the reduced problem satisfied by
$\lambda_r$ is from the projection of the high-fidelity condensed problem. 

\begin{proposition}[Trace estimate]
\label{prop:trace}
Let
$\gamma_k$ be the continuity constant of the sesquilinear form $\mathcal{K}$.  Then
\begin{equation}\label{eq:sum_lambda}
\displaystyle
  \left\| e_\Lambda \right\|_{L^2(\Sigma)}
  \leq\underbrace{\Bigl(1+\frac{\gamma_k}{\beta}\Bigr)
       \inf_{q_r\in \mathcal{Z}_n^\Lambda}\left\| \lambda_h-q_r \right\|_{L^2(\Sigma)}}
       _{\text{best-approx of }\lambda_h\text{ in }\mathcal{Z}_n^\Lambda}
  +\underbrace{\frac{\delta_K\left\| \lambda_h \right\|_{L^2(\Sigma)}+\delta_F}{\beta}}
       _{\text{gap between the two condensations}},
\end{equation}
where $\delta_K, \delta_F>0$ are respectively the operator norms of $\delta \mathbb{K}$ and $\delta \mathbb{F}$.
\end{proposition}
 
\begin{proof}
By  \eqref{eq:commutation}, we get 
\begin{equation}
    \widehat{\mathbb{K}} = \mathbb{Z}_\Lambda^{*}\mathbb{K}\mathbb{Z}_\Lambda-\delta\mathbb{K},
    \qquad
    \widehat{\mathbb{F}} = \mathbb{Z}_\Lambda^{*}\mathbb{F}-\delta\mathbb{F}.
\end{equation}
Fix $q_r\in\mathcal{Z}^\Lambda_n$ and split
$e_\Lambda=(\lambda_h-q_r)+\underbrace{(q_r-\lambda_r)}_{\zeta_r}$; we estimate the
second term. Applying $\widehat{\mathbb{K}}$ to
$\widehat\zeta=\widehat q-\widehat\Lambda$, using $\widehat{\mathbb{K}}\widehat\Lambda=\widehat{\mathbb{F}}$ and $\mathbb{K}\Lambda=\mathbb{F}$, we get
\begin{equation}
   \widehat{\mathbb{K}}\,\widehat\zeta
   =\bigl(\mathbb{Z}_\Lambda^{*}\mathbb{K}\mathbb{Z}_\Lambda-\delta\mathbb{K}\bigr)\widehat q
     -\mathbb{Z}_\Lambda^{*}\mathbb{F}+\delta\mathbb{F}
   =\mathbb{Z}_\Lambda^{*}\mathbb{K}\,(\mathbb{Z}_\Lambda\widehat q-\Lambda)
     -\delta\mathbb{K}\,\widehat q+\delta\mathbb{F}.
\end{equation}
Let $\varphi_r=\bs{Z}_\Lambda\widehat\varphi\in\mathcal{Z}^\Lambda_n$, we get
\begin{equation}
   \widehat\varphi^{*}\widehat{\mathbb{K}}\widehat\zeta
   = \mathcal{K}(q_r-\lambda_h,\varphi_r)
     -\widehat\varphi^{*}\delta\mathbb{K}\,\widehat q
     +\widehat\varphi^{*}\delta\mathbb{F} .
\end{equation}
The first term is bounded by the continuity of $\mathcal{K}$, with constant $\gamma_k$.
For the two others, we introduce
\begin{equation}
   \delta_K:=\sup_{\widehat q,\widehat\varphi\,\in\,\mathbb{C}^{n_\Lambda}\setminus\{0\}}
      \frac{|\widehat\varphi^{\,*}\delta\mathbb{K}\,\widehat q|}
           {\|\widehat q\|\,\|\widehat\varphi\|},
   \qquad
   \delta_F:=\sup_{\widehat\varphi\,\in\,\mathbb{C}^{n_\Lambda}\setminus\{0\}}
      \frac{|\widehat\varphi^{\,*}\delta\mathbb{F}|}{\|\widehat\varphi\|} ,
   \label{eq:deltas}
\end{equation}
the operator norm of $\delta\mathbb{K}$ and the norm of $\delta\mathbb{F}$
on $\mathcal{Z}^\Lambda_n$.  Since $\mathbb{Z}_\Lambda$ has $\mathbb{M}_\Lambda$-orthonormal columns
($\|\widehat q\|=\|q_r\|_{L^2(\Sigma)}$,
$\|\widehat\varphi\|=\|\varphi_r\|_{L^2(\Sigma)}$),
\[
   |\widehat\varphi^{*}\delta\mathbb{K}\,\widehat q|
   \le\delta_K\,\|q_r\|_{L^2(\Sigma)}\,\|\varphi_r\|_{L^2(\Sigma)} ,
   \qquad
   |\widehat\varphi^{*}\delta\mathbb{F}|
   \le\delta_F\,\|\varphi_r\|_{L^2(\Sigma)} .
\]
Altogether,
\begin{equation}
   |\widehat\varphi^{*}\widehat{\mathbb{K}}\widehat\zeta|
   \le\Bigl(\gamma_k\|\lambda_h-q_r\|_{L^2(\Sigma)}
      +\delta_K\|q_r\|_{L^2(\Sigma)}+\delta_F\Bigr)\,
      \|\varphi_r\|_{L^2(\Sigma)} .
\end{equation}
Taking the supremum over $\widehat\varphi\neq0$ and using the inf-sup
constant of $\widehat{\mathbb{K}}$, which is the same $\beta$ as in  \eqref{eq:Kcoerc},
\begin{equation}  
   \beta\,\|\zeta_r\|_{L^2(\Sigma)}
   \le\gamma_k\|\lambda_h-q_r\|_{L^2(\Sigma)}
      +\delta_K\|q_r\|_{L^2(\Sigma)}+\delta_F .
\end{equation}
Choosing $q_r=\Pi_{\mathcal{Z}^\Lambda_n}\lambda_h$, which realizes the infimum
of the first term and satisfies $\|q_r\|_{L^2(\Sigma)}\le\|\lambda_h\|_{L^2(\Sigma)}$,
the triangle inequality $\|e_\Lambda\|_{L^2(\Sigma)}
\le\|\lambda_h-q_r\|_{L^2(\Sigma)}+\|\zeta_r\|_{L^2(\Sigma)}$
gives  \eqref{eq:sum_lambda}.
\end{proof}

\paragraph{3. Primal estimate.} 
Next, we bound the primal error $e_U$. Proposition~\ref{prop:primal} shows that the volume error $e_U$ is controlled by two terms, its best-approximation error in $\mathcal{Z}^U_n$, and the trace
error $e_\Lambda$, propagated through the coupling terms.
 
\begin{proposition}[Primal estimate]
\label{prop:primal}
Let $\gamma_c$ be the continuity constant of $c(q,w):=W^{*}\mathbb{C}Q$, then
\begin{equation}\label{eq:sum_primal}
  \left\| e_U \right\|_{L^2(\Omega)}
  \leq\underbrace{\frac{\gamma_a}{\alpha_a}
  \inf_{v_r\in \mathcal{Z}_n^U}\left\| u_h-v_r \right\|_{L^2(\Omega)}}_{\text{best-approx of }u_h\text{ in }\mathcal{Z}_n^U}
  +\underbrace{C_{\mathrm{coupl}}\left\| e_\Lambda \right\|_{L^2(\Sigma)}}_{\text{propagation of trace error}},
  \qquad
\end{equation}
where $C_{\mathrm{coupl}}=\gamma_c/\alpha_a$.
\end{proposition}
 
\begin{proof}
By  \eqref{eq:local}, the coefficient vector of $u_h$ satisfies $U=\mathbb{A}^{-1}(\mathbb{S}-\mathbb{C}\Lambda)$ while the reduced
back-substitution
$\widehat{\mathbb{A}}\widehat U
=\mathbb{Z}_U^{*}(\mathbb{S}-\mathbb{C}\mathbb{Z}_\Lambda\widehat\Lambda)$ that is $\mathbb{Z}_U \hat{U}=\mathcal{P}(\mathbb{S}-\mathbb{C}\mathbb{Z}_\Lambda\widehat\Lambda)$. Inserting in $E_U$ the coefficient vector of $e_U$, we get 
\[   E_U   =\underbrace{\bigl(\mathbb{A}^{-1}-\mathcal{P}\bigr)      \bigl(\mathbb{S}-\mathbb{C}\Lambda\bigr)}_{\text{(i)}}   +\underbrace{\mathcal{P}\bigl(-\mathbb{C}E_\Lambda\bigr)}_{\text{(ii)}} .\]
Term (i) is  \eqref{eq:cea_local} with $f=\mathbb{S}-\mathbb{C}\Lambda$, it is bounded by$(\gamma_a/\alpha_a)\inf_{v_r}\|u_h-v_r\|_{L^2(\Omega)}$. \\
For term (ii), define the test function $w_r$ with coefficient vector $\mathcal{P}(-\mathbb{C}E_\Lambda)$, it satisfies $a(w_r,v_r)=-c(e_\Lambda,v_r)$ for all $v_r\in\mathcal{Z}^U_n$; choosing $v_r=w_r$ and using coercivity,
\[   \alpha_a\|w_r\|^2_{L^2(\Omega)}\le\operatorname{Re}a(w_r,w_r)   =-\operatorname{Re}c(e_\Lambda,w_r)   \le\gamma_c\|e_\Lambda\|_{L^2(\Sigma)}\|w_r\|_{L^2(\Omega)} .\]
Adding the two contributions gives  \eqref{eq:sum_primal}.
\end{proof}

\paragraph{4. Size of the gaps $\delta_K$ and $\delta_F$.}
Proposition~\ref{prop:gaps} provides upper bounds for the gap terms
$\delta_K$ and $\delta_F$ introduced in \eqref{eq:sum_lambda}. These
bounds show that the gaps are controlled by the distance of the local
volume outputs
\[
\mathbb{A}^{-1}\mathbb{C}\,\mathcal{Z}^\Lambda_n
\qquad\text{and}\qquad
\mathbb{A}^{-1}\mathbb{S}
\]
from the primal reduced space $\mathcal{Z}^U_n$.
\begin{proposition}[Size of the gaps]
\label{prop:gaps}
Let $\gamma_b$ be the continuity constant of the sesquilinear form $b$
associated with $\mathbb{B}$. Then
\begin{equation}
\delta_K
\le
\gamma_b\frac{\gamma_a}{\alpha_a}
\sup_{\|q_r\|_{L^2(\Sigma)}=1}
\operatorname{dist}
\bigl(
\mathbb{A}^{-1}\mathbb{C}q_r,\,
\mathcal{Z}^U_n
\bigr),
\qquad
\delta_F
\le
\gamma_b\frac{\gamma_a}{\alpha_a}
\operatorname{dist}
\bigl(
\mathbb{A}^{-1}\mathbb{S},\,
\mathcal{Z}^U_n
\bigr).
\label{eq:delta_bounds}
\end{equation}
In particular, they both vanish as soon as $\mathbb{A}^{-1}\mathbb{C}\,\mathcal{Z}^\Lambda_n\subset\mathcal{Z}^U_n$ \emph{or} $\mathbb{A}^{-1}\mathbb{S}\in\mathcal{Z}^U_n$, in which case the ROM \eqref{eq:condensed_reduced} and the projected condensed high-fidelity system coincide.  
\end{proposition}

\begin{proof}
The result follows directly from \eqref{eq:cea_local}. For
$q_r\in\mathcal{Z}^\Lambda_n$, applying
\eqref{eq:cea_local} with
$f=\mathbb{C}q_r$ yields the bound for $\delta_K$. Applying the same
estimate with $f=\mathbb{S}$ gives the bound for $\delta_F$.
\end{proof}

\begin{remark}[Primal lifting]
At each training frequency $\sigma_\ell$ for $\ell \in \{1,\ldots,n_{\rm lift}\}$, the two conditions
\[
\mathbb{A}^{-1}(\sigma_\ell)\mathbb{C}\,\mathcal{Z}^\Lambda_n
\subset \mathcal{Z}^U_n,
\qquad
\mathbb{A}^{-1}(\sigma_\ell)\mathbb{S}
\in \mathcal{Z}^U_n
\]
are not independent. Indeed, by construction, the high-fidelity
snapshot $U_h(\sigma_\ell)$ belongs to $\mathcal{Z}^U_n$ and satisfies
\[
U_h(\sigma_\ell)
=
\mathbb{A}^{-1}(\sigma_\ell)\mathbb{S}
-
\mathbb{A}^{-1}(\sigma_\ell)\mathbb{C}\Lambda_h(\sigma_\ell),
\qquad
\Lambda_h(\sigma_\ell)\in\mathcal{Z}^\Lambda_n.
\]
Consequently, if $\mathcal{Z}^U_n$ is enriched with  $
\mathbb{A}^{-1}(\sigma_\ell)\mathbb{C}\,\mathcal{Z}^\Lambda_n
\subset \mathcal{Z}^U_n,
$ then $ \mathbb{A}^{-1}(\sigma_\ell)\mathbb{S}
\in\mathcal{Z}^U_n.$\\
As a result,
both $\delta_K$ and $\delta_F$ vanish at these frequencies. We assess
this lifting strategy numerically in Sections~\ref{sec:marmousi} and
\ref{sec:3d}.
\end{remark}

\begin{remark}[Size of the lifted primal space]
In practice, enforcing the lifting condition
$
\mathbb{A}^{-1}(\sigma_\ell)\mathbb{C}\mathcal{Z}^\Lambda_n
\subset \mathcal{Z}^U_n
$ at $n_{\rm lift}$ selected training frequencies enriches the primal space by at most $n_{\rm lift}n_\Lambda$ additional dimensions, so that $ \dim\mathcal{Z}^U_n \le  n_U+n_{\rm lift}n_\Lambda. $\\
Although these additional basis vectors are inexpensive to compute,
since $\mathbb{A}$ is block diagonal and the required local inverse
matrices are already available from the high-fidelity solves, the
resulting increase in reduced dimension may become significant. This
growth can be controlled in two ways: (i) by performing the lifting
only at $n_{\rm lift}$ selected training frequencies, so that the gaps
vanish exactly at these frequencies and remain small in their
neighborhood, and (ii) by compressing the lifted vectors with a POD
before orthonormalizing them against $\mathcal{Z}^U_n$.

An alternative that keeps the dimension of $\mathcal{Z}^U_n$ fixed is a
Petrov--Galerkin projection of the volume equation, in which the test space
$\mathcal{Y}_n\neq\mathcal{Z}^U_n$ is spanned by the dual liftings
$\mathbb{A}^{-*}(\sigma_\ell)\mathbb{B}^*\mathcal{Z}^\Lambda_n$,
$\ell=1,\dots,n_{\rm lift}$, while the trial space $\mathcal{Z}^U_n$ is left
unchanged. The gaps $\delta_K$ and $\delta_F$ then vanish at the lifting
frequencies without enlarging $\mathcal{Z}^U_n$. However, this approach
requires the analysis of a new inf-sup constant, which is not controlled by
construction since the trial and test spaces are built independently.
Moreover, the reduced system is square only if
$\dim\mathcal{Y}_n=\dim\mathcal{Z}^U_n$  which couples the two reduced spaces. We do not pursue this
approach here and instead consider the two strategies described above.
\end{remark}

\begin{remark}[Relation to localized model order reduction]
The reduced model \eqref{eq:condensed_reduced} is closely related to localized model order reduction methods discussed in the introduction. In both settings, the interior unknowns are approximated in a reduced space and subsequently eliminated by static condensation, resulting in a reduced Schur complement for the interface variables. In the present HDG setting, the volume unknowns play the role of the component interiors, while the trace unknowns on the skeleton play the role of what is referred to as ports in localized model reduction.  In \cite{phuong2013static}, however, the interior reduced space is constructed by lifting the modes of the port reduced space. Consequently, $
\mathbb{A}^{-1}\mathbb{C}\,\mathcal{Z}^\Lambda_n
\subset \mathcal{Z}^U_n $ holds by construction, so that projection and static condensation commute. In the present work, in contrast, the reduced spaces $\mathcal{Z}^U_n$ and $\mathcal{Z}^\Lambda_n$ are constructed independently, and this inclusion is not guaranteed. As a result, condensation and projection do not generally commute, which motivates the analysis presented above.

Moreover, localized methods construct local reduced bases that are valid on a component of the domain and typically rely on local a posteriori error estimators to assess and adapt these bases. For instance, \cite{smetana2015new} develops an estimator accounting for the combined reduction of interior and port variables, while \cite{buhr2017arbilomod,ohlberger2015error} employ localized estimators to drive the adaptive enrichment of local reduced spaces. In the present work, the reduced spaces are instead global over the physical and parameter domains. Accordingly, the estimates established in \cref{prop:trace,prop:primal,prop:gaps} are a priori estimates for the global reduced formulation, rather than localized a posteriori estimators for component-wise enrichment.
\end{remark}

\subsection{Implementation }
\label{sec:overall_algo}

We now summarize the complete offline construction of the reduced
model together with the online evaluation of the obtained ROM.  We combine the ingredients
introduced above: the affine-in-frequency structure, the frequency derivative enrichment
(Section  \ref{sec:basis_construction}), the residual
indicator  \eqref{eq:greedy_indicator}, the weighted POD
compression  \eqref{eq:ROB_extraction}, and the Galerkin
projection  \eqref{eq:reduced_block}. 

The greedy algorithm is initialized with a coarse subset of the
frequency that is sampled given a probability density; at each sampled frequency we perform a
single high-fidelity solve and recover 
its first frequency moments, which are appended to the snapshot set. We perform projection as presented in Section \ref{sec:galerkin} yielding an initial
reduced model that preserves the block structure of the full-order
system. In the greedy stage, the reduced bases are enriched
iteratively: the residual indicator is evaluated over the whole training
set to select the least well-approximated frequency, at which a new
high-fidelity solve and its moments are computed, appended, and
compressed; the loop terminates once the true relative error at
the selected frequency falls below the prescribed tolerance. The complete procedure is
stated in Algorithm  \ref{alg:overall}.

\begin{algorithm}[ht!]
\caption{Offline construction of the frequency ROM at fixed medium}
\label{alg:overall}
\begin{algorithmic}[1]
\State \textbf{Input:} training set $\mathcal{D}^h_\sigma$, initial grid
       $\mathcal{D}^0_\sigma\subset\mathcal{D}^h_\sigma$, moment order
       $m$, tolerance $\varepsilon_{\mathrm{tol}}$
\State \textbf{Output:} reduced bases $(\bs{Z}_U,\bs{Z}_\Lambda)$ and
       reduced operators
       $(\widehat{\mathbb{A}},\widehat{\mathbb{B}},
        \widehat{\mathbb{C}},\widehat{\mathbb{L}},\widehat{\mathbb{S}})$
\Statex
\Statex \textit{Initialization}
\State $\mathcal{S}\gets\varnothing$
       \Comment{snapshot set}
\For{$\sigma\in\mathcal{D}^0_\sigma$}
  \State solve the high-fidelity system  \eqref{eq:condensed} at
         $\sigma$ and recover $X(\sigma)=(U,\Lambda)$
  \State compute the derivatives $\{X^{(j)}\}_{j=0}^{m-1}$
         via  \eqref{eq:taylor_hierarchy}
  \State $\mathcal{S}\gets\mathcal{S}\cup\{X^{(j)}\}_{j=0}^{m-1}$
\EndFor
\State $\bs{Z}_U\gets\mathrm{POD}_{L^2(\Omega)}(\mathcal{S})$,\quad
       $\bs{Z}_\Lambda\gets\mathrm{POD}_{L^2(\Sigma)}(\mathcal{S})$
       \Comment{weighted POD  \eqref{eq:ROB_extraction}}
\State assemble the reduced operators  \eqref{eq:reduced_operators}
       by element-wise projection
\Statex
\Statex \textit{ Greedy enrichment}
\For{$\ell=1,2,\ldots$}
  \State evaluate $\Delta(\sigma)$ over
         $\mathcal{D}^h_\sigma$ \Comment{cost $O(n_U^3)+O(n_\Lambda^3)$ per frequency}
  \State $\sigma^{(\ell)}\gets
         \arg\max_{\sigma\in\mathcal{D}^h_\sigma}\Delta(\sigma)$
    \State \textbf{if} the true relative error at $\sigma^{(\ell)}$ is
         $\le\varepsilon_{\mathrm{tol}}$ \textbf{then break}
  \State solve the high-fidelity system at $\sigma^{(\ell)}$ and
         compute its derivatives $\{X^{(j)}\}_{j=0}^{m-1}$
  \State $\mathcal{S}\gets\mathcal{S}\cup\{X^{(j)}\}_{j=0}^{m-1}$
  \State update $\bs{Z}_U,\bs{Z}_\Lambda$ by weighted POD
         of $\mathcal{S}$ and reassemble the reduced operators
\EndFor 
\Statex
\Statex \textit{ Online evaluation}
\For{$\sigma\in\mathcal{D}_\sigma$}
  \State assemble the system   \eqref{eq:reduced_block} and solve   \eqref{eq:condensed_reduced}

  \State recover $u_r=\bs{Z}_U\widehat{U}$,\quad
         $\lambda_r=\bs{Z}_\Lambda\widehat{\Lambda}$ \Comment{cost $O(n_U^3)+O(n_\Lambda^3)$ per frequency}
\EndFor
\end{algorithmic}
\end{algorithm}

\paragraph{Implementation.} 
Algorithm~\ref{alg:overall} is implemented in the open-source software
\texttt{Hawen} \cite{faucher2021hawen} using a hybrid MPI/OpenMP
parallelization. The implementation exploits the element-wise locality of
the HDG formulation. High-fidelity solves are distributed over MPI
processes, with OpenMP parallelism used for element-wise operations. The construction of the reduced bases, of the reduced operators and of the residual estimator is distributed as follow : each MPI process handles its share of elements, with OpenMP parallelism over them, and the small reduced matrices are summed across the processes.
The main computational steps are summarized below.

\begin{itemize}

\item \emph{High-fidelity solves.}
Thanks to static condensation, only the trace unknown is globally solved,
while the volume unknowns remain local to each element of the mesh. The mesh
is distributed among the MPI processes. For each frequency selected by the
greedy algorithm, on each process, the local matrices $\mathbb{A}_e$ are assembled
and factorized independently for the elements assigned to that process, with
the element-wise operations parallelized using OpenMP. Each process then
assembles its contribution to the condensed global matrix
$\mathbb{K}$ defined in~\eqref{eq:K_def}, which is factorized using the
distributed sparse direct solver MUMPS \cite{amestoy2001fully,amestoy2000mumps}. Once the
trace solution has been obtained, the volume unknowns are recovered
independently on each element through~\eqref{eq:backsubstitution} and parallelized using OpenMP. No global
volume system is therefore assembled or factorized.

\item \emph{Derivative enrichment.}
The computation of the frequency derivatives reuses the factorizations
obtained for the high-fidelity solve. The local contributions to the
right-hand side of~\eqref{eq:moments_condensed} are computed independently
on each element, in parallel with OpenMP threads, and then combined across
MPI processes. The trace derivative
$\Lambda^{(j)}$ is then computed by solving the condensed linear system in~\eqref{eq:moments_condensed} using the factorization of $\mathbb{K}$ obtained for the corresponding high-fidelity solve, while the volume derivatives
$U_e^{(j)}$ are recovered on each element, in parallel with OpenMP threads, using the
factorizations of $\mathbb{A}_e$. Hence, the derivative enrichment requires
no additional factorization.

\item \emph{POD compression.}
The POD compression is  distributed over the MPI processes: each process
keeps its own part of the snapshots and of the reduced basis $\mathbb{Z}_U$,
so that the volume snapshots are never gathered on a single process.

\item \emph{Local projection of the reduced operators.}
Let $\mathbb{Z}_{U,e}$ denote the restriction of $\mathbb{Z}_U$ to the element
$K_e$, and let $\mathcal{R}_e\mathbb{Z}_\Lambda$ denote the restriction of the
trace basis to the faces of $K_e$. The reduced operators are assembled
directly from element-wise contributions. For instance,
\begin{equation}
\widehat{\mathbb{A}}_k
=
\sum_e
\mathbb{Z}_{U,e}^{*}\,\mathbb{A}_{k,e}\,\mathbb{Z}_{U,e},
\quad k=0,1,
\qquad
\widehat{\mathbb{C}}
=
\sum_e
\mathbb{Z}_{U,e}^{*}\,\mathbb{C}_e\,\mathcal{R}_e\mathbb{Z}_\Lambda,
\label{eq:local_projection}
\end{equation}
with analogous expressions for
$\widehat{\mathbb{B}}$, $\widehat{\mathbb{L}}$, and
$\widehat{\mathbb{S}}$. Thus, the contribution associated with an element
depends only on its local matrices and on the restrictions of the reduced
bases to that element. These contributions can be computed independently
and accumulated into the reduced matrices. In the implementation, the
cell loop is parallelized with OpenMP on each MPI process, and the reduced matrices are then summed across the processes.

Importantly, the global uncondensed operator~\eqref{eq:full_block} is never
assembled. The reduced operators are constructed directly from the local
HDG matrices, as described in Remark~\ref{rem:local_proj}. This preserves
the element-wise structure of the HDG discretization and avoids the memory
and computational cost associated with the assembly of the global volume
operator. Moreover, $\widehat{\mathbb{A}}_0$ and
$\widehat{\mathbb{A}}_1$ are stored separately so that the affine
dependence
\[
\widehat{\mathbb{A}}(\sigma)
=
\widehat{\mathbb{A}}_0+\sigma\widehat{\mathbb{A}}_1
\]
is preserved explicitly.

\item \emph{Residual indicator.}
The quantities entering the residual indicator $\Delta$ defined in equation \cref{eq:greedy_indicator} are constructed in the same
element-wise manner. In particular, the products
\[
\mathbb{A}_0\mathbb{Z}_U,\qquad
\mathbb{A}_1\mathbb{Z}_U,\qquad
\mathbb{C}\mathbb{Z}_\Lambda
\]
are evaluated element by element and using OpenMP, on each MPI process. Since the volume mass matrix
$\mathbb{M}_U$ is block diagonal, the corresponding volume contributions
to the Gram matrices~\eqref{eq:gram_volume} can be accumulated locally and summed across the processes. Similarly, the latter holds for the products involving $\mathbb{B}$ and $\mathbb{L}$ and the norm matrix
$\mathbb{M}_\Lambda^{-1}$.

The resulting Gram matrices involve only the reduced order operators and are
therefore small dense matrices. They are recomputed after each update of the reduced
bases. Once these matrices have been constructed, the
evaluation of the residual indicator over the training set
$\mathcal{D}^h_\sigma$ requires only reduced order computations: for each
candidate frequency, the reduced problem is solved and the residual norm
is evaluated through the quadratic form~\eqref{eq:res_quadratic}. 
\end{itemize}

Overall, the high-fidelity solves and derivative computations rely on
MPI-distributed sparse linear algebra, combined with OpenMP parallelism for
element-wise operations. The construction of the reduced operators and
the residual estimator exploits the same element-wise structure and is
parallelized over the MPI processes and over the cells using OpenMP. In particular, the global
uncondensed volume operator ~\eqref{eq:full_block}  is never assembled during the reduced model
construction. This implementation significantly reduces the offline
computational cost in the two- and three-dimensional experiments reported
in Sections~\ref{sec:marmousi} and~\ref{sec:3d}.

\section{Two-dimensional model: Marmousi}
\label{sec:marmousi}
We validate the frequency-domain ROM construction presented in
Section~\ref{sec:ROM} on a fixed heterogeneous medium $\mathbf{m}$, using the
two-dimensional Marmousi model. We consider the acoustic wave propagation
problem \eqref{eq:acoustic_strong} on a rectangular domain $\Omega$, with a
Neumann boundary condition imposed on the top boundary and absorbing boundary
conditions on the remaining three sides. The problem is discretized using the
HDG scheme described in Section~\ref{sec:discretization}, with polynomial
degrees $p=5$. The Marmousi model, provided by IFPEN,
is considered with constant density $\rho=1$.

\begin{figure}[ht!]
    \centering
    \subfloat[2D Marmousi wave speed model of size 9.2$\times$3 km$^2$.\label{fig:marmousi}]
    {\includegraphics[scale=1.2]{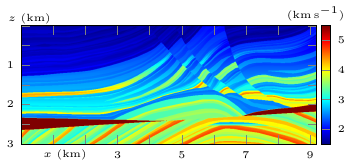}}

    \subfloat[
    Real part of the pressure field at 
    frequency $\omega=2\pi\times 8.51$
    (left)
    and $2\pi\times 9.96$
    (right)
    using the high-fidelity solver, for a source located at $(1.5,0.1)$ km.
    ]{
    \includegraphics[scale=1]{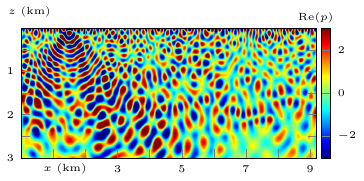} \hspace*{3em}
    \includegraphics[scale=1]{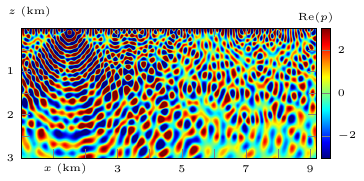}
    }

    \caption{Two-dimensional Marmousi problem using a frequency interval from 8.51 to 9.96 Hz.}
    \label{fig:pressure_mamrousi}
\end{figure}

The frequency band of interest is
$\mathcal{D}_\omega=[8.5,10]\times 2\pi$,
with damping parameter $s=0.001$. The reduced model is constructed using
Algorithm~\ref{alg:overall} on a computational mesh comprising
$N_e=15{,}000$ elements leading to $N_U=992,376$ and $N_\Lambda=142,776$. A sharp localized Gaussian source is located at
$\boldsymbol{x}=(1500,100)$, near the top-left corner of the domain. Figure~\ref{fig:pressure_mamrousi}
illustrates the resulting high-fidelity pressure field $p$ at frequencies
$\omega=8.51\times 2\pi$ and $\omega=9.96\times 2\pi$, highlighting the influence of the heterogeneous
medium and the frequency on wave propagation.

We conduct an extensive numerical investigation of the main ingredients of
Algorithm~\ref{alg:overall}, together with the error analysis presented in
Section~\ref{sec:error_bounds}, before assessing the overall performance of the
resulting reduced model~\eqref{eq:condensed_reduced}. The investigation is
organized as follows:
\begin{itemize}
    \item First, we investigate the benefits and computational cost of enriching
    the reduced basis with Fréchet derivatives. To this end, we construct
    reduced models using plain POD, both with and without sensitivity
    enrichment, and assess the resulting approximation accuracy and
    computational cost.

    \item Second, we assess the accuracy of the error bounds derived in
    Section~\ref{sec:error_bounds}. In particular, we investigate the lifting
    strategy introduced in that section on a POD and compare it to a POD
    without lifting, in order to determine whether lifting effectively
    decouples the approximation errors associated with the two unknowns.

    \item Third, we evaluate the performance of the residual-based error indicator
    \eqref{eq:greedy_indicator}. Specifically, we use the reduced model
    constructed with plain POD to evaluate the cost of building  the indicator and compare its
    frequency dependence with that of the true approximation error. This
    allows us to assess whether the indicator identifies the same critical
    frequencies as the actual error.

    \item Finally, having validated the preceding ingredients, we combine them
    within Algorithm~\ref{alg:overall} to construct the reduced model and
    assess its overall performance:
    \begin{itemize}
        \item the convergence of the greedy algorithm in terms of in-sample
        error decay, number of iterations, and computational cost;

        \item the offline and online computational costs, including the cost
        of each component of the greedy procedure
        Algorithm~\ref{alg:overall}, as well as the cost of evaluating
        out-of-sample parameters and the resulting speedup per query;

        \item the accuracy of the reduced model on out-of-sample testing
        parameters, using the error metrics defined below.
    \end{itemize}
\end{itemize}

\paragraph{Computational setup.}
All 2D Marmousi computations were performed on a single node of  PlaFRIM platform \cite{PlaFRIM}, with two 24-core AMD Zen4 processors and 256\,GB of memory. We use a hybrid MPI/OpenMP parallelization (Cf. \cref{sec:overall_algo}) with $4$ MPI processes and $12$ OpenMP threads.

\paragraph{Error metrics.}

For a given frequency $\omega$, we assess the reduced model against the
high-fidelity solution using relative $L^2$ projection and approximation errors
for each unknown. The projection error measures the best approximation
achievable within the corresponding reduced space, independently of the
reduced order solve. It therefore isolates the approximation properties of the
reduced basis itself. For the volume unknown $u_h$, it is defined as
\begin{equation}
    E^U_{\mathrm{proj}}(\omega)
    = \frac{\|u_h(\omega) - \Pi_{\mathcal{Z}^U_n} u_h(\omega)\|_{L^2(\Omega)}}
           {\|u_h(\omega)\|_{L^2(\Omega)}},
    \label{eq:proj_ern_U}
\end{equation}
and, analogously, for the trace unknown $\lambda_h$,
\begin{equation}
    E^\Lambda_{\mathrm{proj}}(\omega)
    = \frac{\|\lambda_h(\omega) - \Pi_{\mathcal{Z}^{\Lambda}_n} \lambda_h(\omega)\|_{L^2(\Sigma)}}
           {\|\lambda_h(\omega)\|_{L^2(\Sigma)}}.
    \label{eq:proj_err_L}
\end{equation}
Here, $\mathcal{Z}^{\bullet}_n$ denotes the reduced space associated with the
corresponding unknown, while
$\Pi_{\mathcal{Z}^U_n}u_h(\omega)$ and
$\Pi_{\mathcal{Z}^{\Lambda}_n}\lambda_h(\omega)$ denote the respective
$L^2$-orthogonal projections onto these spaces.

In contrast, the approximation error measures the actual accuracy of the
reduced solution obtained by solving the reduced problem. For the volume
unknown $u_r(\omega)\in\mathcal{Z}^U_n$ and the trace unknown
$\lambda_r(\omega)\in\mathcal{Z}^{\Lambda}_n$, the relative approximation
errors are defined as
\begin{equation}
    E^U_{\mathrm{rel}}(\omega)
    = \frac{\|u_h(\omega) - u_r(\omega)\|_{L^2(\Omega)}}
           {\|u_h(\omega)\|_{L^2(\Omega)}},
    \label{eq:rel_ern_U}
\end{equation}
and
\begin{equation}
    E^\Lambda_{\mathrm{rel}}(\omega)
    = \frac{\|\lambda_h(\omega) - \lambda_r(\omega)\|_{L^2(\Sigma)}}
           {\|\lambda_h(\omega)\|_{L^2(\Sigma)}}.
    \label{eq:rel_err_L}
\end{equation}

Comparing the projection and approximation errors quantifies the additional
error introduced by the reduced order solve beyond that associated with the
finite-dimensional reduced spaces. We therefore introduce the effectivity
index
\begin{equation}
    \theta^U_{\mathrm{eff}}(\omega)
    = \frac{E^U_{\mathrm{rel}}(\omega)}{E^U_{\mathrm{proj}}(\omega)}
    \;\geq 1,
    \label{eq:eff_index}
\end{equation}
and, analogously, $\theta^\Lambda_{\mathrm{eff}}(\omega)$ for the trace
unknown. A value close to $1$ indicates that the reduced solution is
quasi-optimal in the corresponding reduced space, meaning that the reduced
solve introduces a small additional error relative the projection error.

\subsection{Enhancement of the ROB via frequency derivatives}

We assess the benefit of enriching the training snapshots with their first order
frequency derivatives , at a matched offline
training budget of $n_{\mathrm{train}}=30$ frequency points uniformly
distributed over $\mathcal{D}_\omega=[8.5,10] \times 2\pi\,\mathrm{ }$, and out-of-sample test set of
$n_{\mathrm{test}}=10$ frequencies.

At the same POD tolerance ($10^{-10}$), enrichment increases the retained basis size from
$n_U=n_\Lambda=30$ to $n_U=n_\Lambda=57$, indicating that the
derivative information is not redundant with the primal snapshots and
meaningfully enriches the reduced space. Table~\ref{tab:enrichment}
reports the relative projection error $E_{\mathrm{proj}}$, the relative $L^2$
error $E_{\mathrm{rel}}$, and the effectivity index $\theta_{\mathrm{eff}}$
for both variables, with and without enrichment.

We observe an improvement of two orders of magnitude in the mean projection
error for both variables $u_h$ and $\lambda_h$. For example, the mean
projection error for $u_h$ decreases from
$1.21\times10^{-2}$ to $1.70\times10^{-4}$, demonstrating that the
derivative snapshots provide additional information beyond that contained
in the primal snapshots. A comparable improvement is observed in the ROM
error: the maximum error in $u_h$ decreases from
$5.15\times10^{-2}$ to $9.83\times10^{-4}$, with $\lambda_h$ exhibiting
the same behavior. 

This improvement is obtained at a moderate additional offline cost.
Computing the $30$ high-fidelity snapshots together with their frequency
derivatives requires $24.08\,\mathrm{s}$ in total, of which only
$1.97\,\mathrm{s}$, corresponding to  $8.2\%$ of the
total training cost, is associated with the derivative computation. This
limited additional cost results from reusing the same
factorizations for the primal and differentiated systems as described in \cref{sec:basis_construction}.

The effectivity index is close to the same range with and without
enrichment, with $\theta_{\mathrm{eff}} \approx 1.6$--$2.2$ in both cases.
This indicates that the reduced solution remains close to the best
approximation achievable within the corresponding reduced space. In
particular, the moment-enrichment does not 
deteriorate  the quasi-optimality of the
Galerkin reduced solution. This justifies the use of the Galerkin
formulation adopted here without introducing a Petrov--Galerkin
projection. As a conclusion, we retain the frequency-sensitivity enrichment for
the remainder of the numerical investigation.

\begin{table}[H]
\centering
\small
\begin{tabular}{lcccccc}
\toprule  & rank $n$ 
 & mean $E_{\mathrm{proj}}$ & max $E_{\mathrm{proj}}$ & mean $E_{\mathrm{rel}}$ & max $E_{\mathrm{rel}}$ & $\theta_{\mathrm{eff}}$ (mean/max) \\
\midrule
\multicolumn{7}{l}{\textit{Variable $u_h$}} \\
without freq. sensitivity & $30$  & $1.21\times10^{-2}$ & $2.97\times10^{-2}$ & $2.10\times10^{-2}$ & $5.15\times10^{-2}$ & $1.74$ / $1.99$ \\
with freq. sensitivity & $57$   & $1.70\times10^{-4}$ & $6.24\times10^{-4}$ & $2.65\times10^{-4}$ & $9.83\times10^{-4}$ & $1.55$ / $2.23$ \\
\addlinespace
\multicolumn{7}{l}{\textit{Variable $\lambda_h$}} \\
without freq. sensitivity & $30$  & $1.21\times10^{-2}$ & $2.97\times10^{-2}$ & $2.10\times10^{-2}$ & $5.14\times10^{-2}$ & $1.74$ / $1.99$ \\
with freq. sensitivity    & $57$  & $1.70\times10^{-4}$ & $6.23\times10^{-4}$ & $2.65\times10^{-4}$ & $9.83\times10^{-4}$ & $1.55$ / $2.24$ \\
\bottomrule
\end{tabular}
\caption{Enhancement of the approximation properties using frequency derivatives.}
\label{tab:enrichment}
\end{table}

\subsection{Analysis of error bounds}
We now assess the accuracy of the bounds derived in \cref{sec:error_bounds}. We first examine the gap terms appearing in Proposition~\ref{prop:trace}, whose bounds are given by~\eqref{eq:delta_bounds}. The trace basis $\boldsymbol{Z}_\Lambda$ is kept fixed throughout the experiment, with a POD tolerance of $10^{-10}$, resulting in $E^\Lambda_{\rm proj}=4.4126\times10^{-5}$. The volume basis $\boldsymbol{Z}_U$ is varied by changing its POD tolerance ${\rm tol}_U$. For each tolerance, the ROM is constructed both without lifting and with the lifting condition
$\mathbb{A}^{-1}(\sigma_\ell)\mathbb{C}\,\mathcal{Z}^\Lambda_n\subset\mathcal{Z}^U_n$
enforced at $n_{\rm lift}=3$ frequencies $\sigma_\ell$. The resulting gap terms, projection errors, and relative errors are reported in Table~\ref{tab:gaps}.

\begin{table}[H]
\centering
\small
\setlength{\tabcolsep}{4pt}
\begin{tabular}{lccccccc}
\toprule
 ${\rm tol}_U$ & $\delta_K$ & $\delta_F$
 & $E^{U}_{\mathrm{proj}}$ & $E^{\Lambda}_{\mathrm{proj}}$
 & $E^{U}_{\mathrm{rel}}$ & $E^{\Lambda}_{\mathrm{rel}}$
 & $\theta^{U}_{\mathrm{eff}}$ \\
\midrule
\multicolumn{8}{l}{\textit{without lifting}} \\
$10^{-6}$ & $1.00$ & $3.19\times10^{-2}$ & $1.4226\times10^{-3}$ & $4.4126\times10^{-5}$ & $3.1895\times10^{-3}$ & $3.1894\times10^{-3}$ & $2.2$ \\
$10^{-8}$ & $1.00$ & $3.00\times10^{-2}$ & $2.5073\times10^{-4}$ & $4.4126\times10^{-5}$ & $5.4316\times10^{-4}$ & $5.4321\times10^{-4}$ & $2.2$ \\
$10^{-10}$ & $8.74\times10^{-5}$ & $2.86\times10^{-2}$ & $4.4143\times10^{-5}$ & $4.4126\times10^{-5}$ & $1.2692\times10^{-4}$ & $1.2688\times10^{-4}$ & $2.9$ \\
$10^{-12}$ & $5.80\times10^{-5}$ & $2.86\times10^{-2}$ & $6.5862\times10^{-6}$ & $4.4126\times10^{-5}$ & $1.0149\times10^{-4}$ & $1.0146\times10^{-4}$ & $15.4$ \\
$10^{-14}$ & $5.80\times10^{-5}$ & $2.86\times10^{-2}$ & $1.6414\times10^{-6}$ & $4.4126\times10^{-5}$ & $1.0149\times10^{-4}$ & $1.0146\times10^{-4}$ & $61.8$ \\
$10^{-16}$ & $5.80\times10^{-5}$ & $2.86\times10^{-2}$ & $9.8328\times10^{-7}$ & $4.4126\times10^{-5}$ & $1.0149\times10^{-4}$ & $1.0146\times10^{-4}$ & $103.2$ \\
\addlinespace
\multicolumn{8}{l}{\textit{with lifting}, $n_{\rm lift}=3$} \\
$10^{-6}$ & $1.58\times10^{-13}$ & $6.23\times10^{-11}$ & $3.5083\times10^{-5}$ & $4.4126\times10^{-5}$ & $1.0070\times10^{-4}$ & $1.0064\times10^{-4}$ & $2.9$ \\
$10^{-8}$ & $2.38\times10^{-13}$ & $6.21\times10^{-11}$ & $1.6967\times10^{-5}$ & $4.4126\times10^{-5}$ & $1.0070\times10^{-4}$ & $1.0064\times10^{-4}$ & $5.9$ \\
$10^{-10}$ & $1.80\times10^{-13}$ & $6.18\times10^{-11}$ & $9.0753\times10^{-6}$ & $4.4126\times10^{-5}$ & $1.0070\times10^{-4}$ & $1.0064\times10^{-4}$ & $11.1$ \\
$10^{-12}$ & $2.06\times10^{-13}$ & $6.12\times10^{-11}$ & $4.6320\times10^{-6}$ & $4.4126\times10^{-5}$ & $1.0070\times10^{-4}$ & $1.0064\times10^{-4}$ & $21.7$ \\
$10^{-14}$ & $1.80\times10^{-13}$ & $6.10\times10^{-11}$ & $1.5667\times10^{-6}$ & $4.4126\times10^{-5}$ & $1.0070\times10^{-4}$ & $1.0064\times10^{-4}$ & $64.3$ \\
$10^{-16}$ & $1.81\times10^{-13}$ & $6.10\times10^{-11}$ & $9.4277\times10^{-7}$ & $4.4126\times10^{-5}$ & $1.0070\times10^{-4}$ & $1.0064\times10^{-4}$ & $106.8$ \\
\bottomrule
\end{tabular}
\caption{Gaps $\delta_K$, $\delta_F$, projection errors and relative errors, with and without lifting, at frozen $\boldsymbol{Z}_\Lambda$.}
\label{tab:gaps}
\end{table}

The lifted configurations  confirm the decoupling predicted by Proposition~\ref{prop:trace}. Enforcing the lifting condition reduces both gap terms to numerical precision, with $\delta_K\approx10^{-13}$ and $\delta_F\approx10^{-10}$ throughout the experiment. Consequently, the trace error becomes independent of the volume reduced space: $E^\Lambda_{\rm rel}$ remains equal to $1.01\times10^{-4}$ for all six lifted configurations, despite the variation of $E^U_{\rm proj}$ by more than one order of magnitude. This contrasts with the unlifted configurations, for which the trace error varies with the volume basis through the gap terms. This behavior is consistent with \cref{eq:sum_lambda}: when $\delta\mathbb{K}$ and $\delta\mathbb{F}$ vanish, the reduced trace problem~\eqref{eq:condensed_reduced} corresponds to the projection of the high-fidelity condensed problem onto $\mathcal{Z}^\Lambda_n$ and is therefore independent of $\mathcal{Z}^U_n$. The trace error is consequently governed only by the approximation properties of $\boldsymbol{Z}_\Lambda$ giving $\theta^\Lambda_{\rm eff}=2.28$.

The lifted configurations also illustrate the behavior predicted by Proposition~\ref{prop:primal}. The relative volume error remains constant at $E^U_{\rm rel}=1.0070\times10^{-4}$ across all volume-basis tolerances, while the volume projection error decreases by more than one order of magnitude. Thus, once the volume best-approximation error becomes smaller than the contribution propagated from the trace variable, further enrichment of $\mathcal{Z}^U_n$ does not improve the total volume error. The latter is therefore controlled by the trace approximation and its propagation to the volume variable, in agreement with the structure of the bound~\eqref{eq:sum_primal}.

The unlifted configurations provide an additional observation. As the volume basis is enriched, the gap terms decrease and the errors become comparable to those obtained with lifting. Since the volume snapshots satisfy
$u_h(\omega_\ell)=\mathbb{A}^{-1}\mathbb{S}-\mathbb{A}^{-1}\mathbb{C}\lambda_h(\omega_\ell)$, so that, up to the source term, they contain the lifting of the corresponding trace snapshots. Consequently, a POD basis constructed from the volume snapshots naturally captures a significant portion of $\mathbb{A}^{-1}\mathbb{C}\,\mathcal{Z}^\Lambda_n$. For our test case, this makes the enforcement of the lifting condition unnecessary once a sufficiently rich volume basis is used.

For this reason, we retain Algorithm~\ref{alg:overall} as it is without explicitly enforcing the lifting condition in the remainder of the numerical experiments, where relatively tight POD tolerances are employed. However, the lifting strategy remains relevant in settings where the volume basis is constructed from limited or insufficient snapshot data, or for saddle-point problems in which the volume and trace variables have distinct structures and the volume snapshots do not naturally capture the lifting of the trace space. In such cases, the gap terms $\delta_K$ and $\delta_F$ can provide a significant contribution to the ROM error, as illustrated by the coarse unlifted configurations in Table~\ref{tab:gaps}.

Overall, the lifting strategy effectively decouples the trace approximation from the volume reduced space, as predicted by Proposition~\ref{prop:trace}, and fixes the trace error at the best-approximation level determined by $\boldsymbol{Z}_\Lambda$. Second, for our model problem, this decoupling is not necessary because the structure of the snapshots already makes the gap terms small enough for sufficiently rich volume spaces. Finally,  as predicted by \cref{prop:primal}, the volume error is bounded by both its best-approximation term and the propagated trace error.

\subsection{Behavior of the error indicator}

\begin{figure}[ht!]
\centering
\includegraphics[scale=.35]{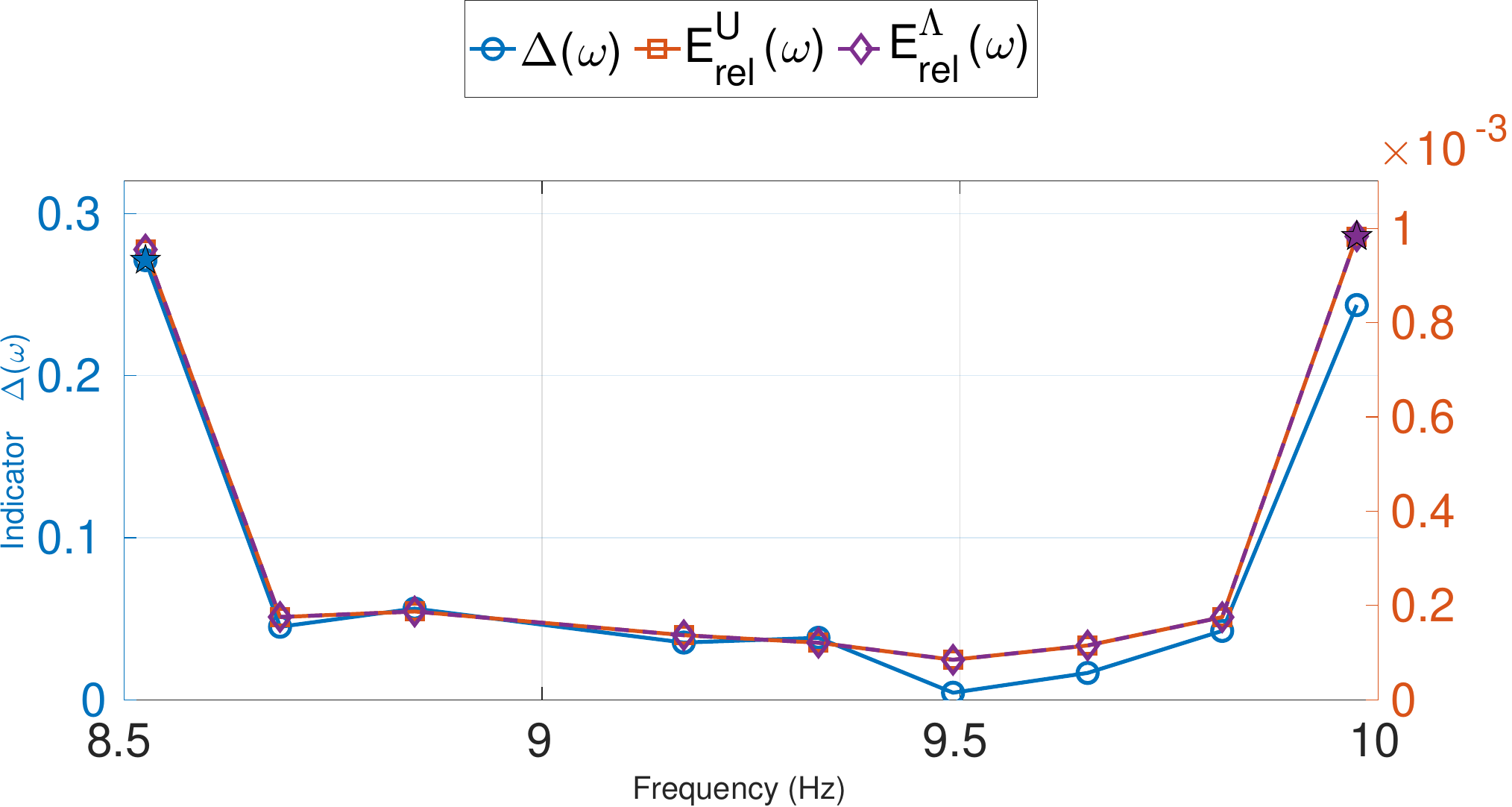}
\caption{Indicator $\Delta(\omega)$ and true relative error
$E^U_{\mathrm{rel}}(\omega)$ on the test set; the marked points indicate $\arg\max_\omega
\Delta(\omega)$ and $\arg\max_\omega E^U_{\mathrm{rel}}(\omega)$.}
\label{fig:indicator_test1}
\end{figure}

In this section, we evaluate the error indicator $\Delta(\omega)$ defined in \eqref{eq:greedy_indicator}. On a plain POD run with a uniform training grid of
$n_{\mathrm{train}}=30$ frequencies, we compare $\Delta(\omega)$ with the relative error $E_{\mathrm{rel}}(\omega)$ defined in \cref{eq:rel_ern_U}--\cref{eq:rel_err_L}, on a test set of $n_{\mathrm{test}}=10$ frequencies. Figure~\ref{fig:indicator_test1} shows that the residual-based indicator closely follows the behavior of the relative error. Both quantities identify the same pair of critical frequencies, $\omega=8.53$ and $\omega=9.97$, but rank them in opposite order: the indicator is maximal at $\omega=8.53$, where it reaches $2.71\times10^{-1}$, whereas the true error is maximal at $\omega=9.97$. Although this leads to a suboptimal selection in the greedy step, the frequency selected by the indicator still captures $97\%$ of the maximal error.

We also observe that the relative errors of the two variables $u_h$ and $\lambda_h$ coincide. This is expected since the skeleton variable $\lambda_h$ represents the numerical trace of the pressure component of $u_h$. The indicator and the relative error, however, are not expected to coincide in magnitude. First, the theoretical bound is not necessarily sharp and  involves a scaling with the inf-sup stability constant. Second, the bound is in terms of the absolute error, whereas the quantity reported here is a relative error. The indicator as it is defined in \cref{eq:greedy_indicator} should be interpreted as a measure for identifying the frequencies where the ROM error is largest, rather than as a sharp estimate of the relative error magnitude.

\subsection{Convergence of the greedy algorithm}

The previous sections validated separately the steps involved in the basis and ROM construction described in \cref{sec:ROM}. We assess now the overall performance of \cref{alg:overall}. Starting from $n_{\mathrm{init}}=30$ uniformly distributed snapshots and a candidate set of $n_{\mathrm{cand}}=100$ frequencies, we set the tolerance of the POD compression step in the algorithm to $10^{-12}$ for both variables. The greedy loop is terminated when the relative error falls below $5\times10^{-4}$.

We report the evolution of the indicator $\Delta(\omega)$ defined in \cref{eq:greedy_indicator} and the true relative error $E_{\mathrm{rel}}(\omega)$ at each iteration of the greedy run in Algorithm \ref{alg:overall}. This allows us to assess whether reducing the indicator $\Delta(\omega)$ along the greedy loop effectively reduces the true error, thereby evaluating the quasi-optimality of the greedy selection strategy. The loop terminates after three iterations.

\begin{figure}[ht!]
\centering
\includegraphics[width=0.4\linewidth]{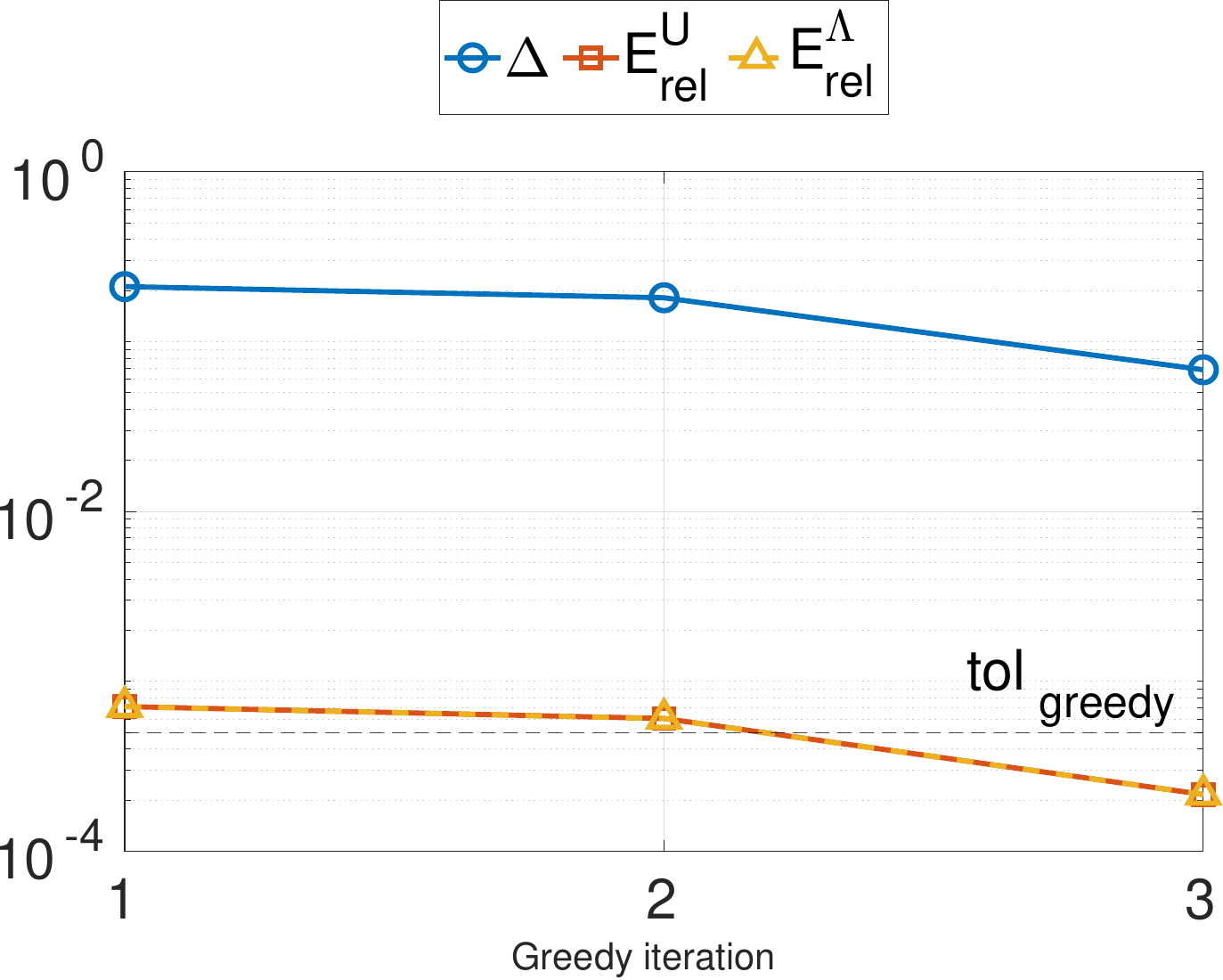}
\caption{Evolution of the indicator $\Delta$ and the relative errors
$E^U_{\mathrm{rel}}$, $E^\Lambda_{\mathrm{rel}}$ along the greedy
iterations, the loop stops when $E^U_{\rm rel}\, , \,  E^\Lambda_{\rm rel} \le \varepsilon_{\mathrm{tol}}=5\times10^{-4}$.}
\label{fig:indicator_test2}
\end{figure}

Figure \ref{fig:indicator_test2} shows that both $\Delta$ and $E_{\mathrm{rel}}$ decrease monotonically, with the relative error decreasing from $7.11\times10^{-4}$ to $2.16\times10^{-4}$ at the selected greedy points. Evaluating the indicator over the $100$ candidates costs $3\times10^{-3}\mathrm{s}$ per iteration, which represents a negligible fraction of the total iteration cost and makes a dense exploration of the parameter domain affordable.

\subsection{Offline and online costs}

Offline, constructing the ROM requires $33$ high-fidelity (HF) solves
($n_{\mathrm{init}}=30$ initial snapshots plus $3$ greedy iterations),
resulting in $26.02\,\mathrm{s}$ of HF solve time ($24.10\,\mathrm{s}$
for the initial training set and $1.92\,\mathrm{s}$ within the loop), out of a
total offline cost of $45.63\,\mathrm{s}$. The remaining
$19.61\,\mathrm{s}$ is spent on POD compression, Gram matrix
assembly, projection, and indicator evaluation across the greedy
loop (Cf. \cref{sec:basis_construction}). The offline cost is therefore dominated by the HF solves,
while the compression and projection procedures account for about
$43\%$ of the total offline cost.

Online, evaluating the ROM at a new frequency costs on average
$3.3\times10^{-4}\,\mathrm{s}$, compared with $0.637\,\mathrm{s}$ for a
full-order solve, yielding a speedup of approximately
$1952\times$ (Table~\ref{tab:speedup}). The offline training is
therefore recovered after approximately $72$ online queries.

\begin{table}[htbp]
\centering
\begin{tabular}{lc}
\toprule
Offline: number of HF solves & $33$ \\
Offline: HF solve time & $26.02$ s \\
Offline: POD + projection + residual indicator build & $19.61$ s \\
Offline: total time & $45.63$ s \\
\midrule
Online: FOM time per query & $0.637$ s \\
Online: ROM time per query & $3.3\times10^{-4}$ s \\
Online: speedup per query & $1952\times$ \\
\bottomrule
\end{tabular}
\caption{Offline and online computational costs of the final greedy ROM
($3$ iterations, $n_U=n_\Lambda=65$), evaluated on
$n_{\mathrm{test}}=50$ independent test frequencies.}
\label{tab:speedup}
\end{table}

\subsection{Accuracy: FOM vs ROM}

The final ROM, obtained after $3$ greedy iterations with
$n_U=n_\Lambda=65$, is validated on an independent test grid of
$n_{\mathrm{test}}=50$ frequencies uniformly sampled over
$\mathcal{D}_\omega=[8.5,10] \times 2\pi$, disjoint from the training
frequency set. Table \ref{tab:accuracy} reports the projection errors,
with a mean value of $2.45\times10^{-5}$, indicating that the greedy
procedure effectively selects snapshots that provide an accurate
reduced basis. This demonstrates the ability of the  basis to
represent the solution over the test set, for solutions with
oscillatory behavior and heterogeneity of the medium. It therefore shows the effectiveness of  combining the greedy selection with sensitivity
enrichment strategies.

Table \ref{tab:accuracy} also reports the relative $L^2$ errors over
the test set, separately for $(p_h,\boldsymbol v_h)$ and $\lambda_h$. For
$p_h$, the mean (resp. maximum) relative error is
$4.41\times10^{-5}$ (resp. $1.66\times10^{-4}$), and the velocity
components $v_{h,x}$ and $v_{h,z}$ have errors of the same order of
magnitude. The corresponding error values for $\lambda_h$ are identical to
those obtained for $p_h$. In all cases, the maximum relative errors
remain below the greedy tolerance
$\varepsilon_{\mathrm{tol}}=5\times10^{-4}$ used to terminate the
training. Finally, the ratio between the relative and projection
errors remains close to $1.8$ for all fields, indicating that the
Galerkin ROM has an approximation quality close to that of the
corresponding reduced bases.

\begin{table}[H]
\centering
\begin{tabular}{lcccc}
\toprule
 & mean $E_{\mathrm{proj}}$ & max $E_{\mathrm{proj}}$ & mean $E_{\mathrm{rel}}$ & max $E_{\mathrm{rel}}$ \\
\midrule
$p_h$      & $2.45\times10^{-5}$ & $8.42\times10^{-5}$ & $4.41\times10^{-5}$ & $1.66\times10^{-4}$ \\
$v_{h,x}$  & $2.46\times10^{-5}$ & $8.83\times10^{-5}$ & $4.28\times10^{-5}$ & $1.68\times10^{-4}$ \\
$v_{h,z}$  & $2.47\times10^{-5}$ & $8.32\times10^{-5}$ & $4.31\times10^{-5}$ & $1.56\times10^{-4}$ \\
$\lambda_h$  & $2.45\times10^{-5}$ & $8.42\times10^{-5}$ & $4.41\times10^{-5}$ & $1.66\times10^{-4}$ \\
\bottomrule
\end{tabular}
\caption{Projection and relative $L^2$ errors of the final greedy ROM  on an independent test grid of $n_{\mathrm{test}}=50$ frequencies over $\mathcal{D}_\omega=[8.5,10]\times 2 \pi$  .}
\label{tab:accuracy}
\end{table}

\begin{figure}[ht!]
    \centering
    \includegraphics[scale=1.2]{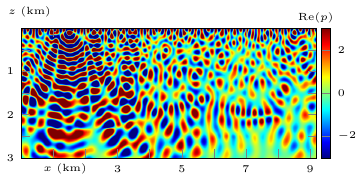} \hspace*{2em}
    \includegraphics[scale=1.2]{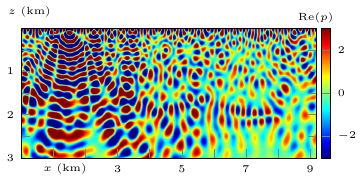}

    \includegraphics[scale=1.25]{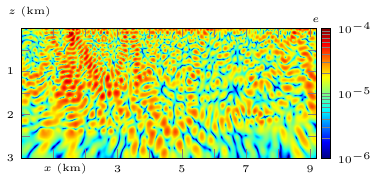}
    \caption{Pressure wavefield  at the worst test frequency $\omega=8.69 \times 2\pi$: (top left) high-fidelity solution, (top right) reduced solution and (bottom) relative error between the two in $\log_{10}$ scale.}
    \label{fig:hf_vs_rom_marmousi}
\end{figure}

Figure~\ref{fig:hf_vs_rom_marmousi} compares the high-fidelity and reduced pressure field at the test frequency of largest error, $\omega=8.69\times2\pi$. We note that the fields are visually identical and that the associated relative error is distributed throughout the computational domain, with values ranging from $10^{-6}$ to $10^{-4}$.

\subsection{Influence of the frequency band}

We repeat the greedy construction on two other frequency bands,
$\mathcal{D}_\omega=[5,7]\times2\pi$ and $\mathcal{D}_\omega=[13,15]\times2\pi$,
with the same settings as above. Each reduced model is then evaluated on $n_{\mathrm{test}}=50$ independent
test frequencies of its band. Table~\ref{tab:bands} reports, for each band,
the number of HF solves, the size of the reduced bases, the maximum relative
errors over the test set, the offline cost and the online speedup.

\begin{table}[H]
\centering
\begin{tabular}{lccccccc}
\toprule
band (Hz) & HF solves & iterations & $n_U=n_\Lambda$ & max $E^U_{\mathrm{rel}}$ & max $E^\Lambda_{\mathrm{rel}}$ & offline time & speedup \\
\midrule
$[5,7]$    & $42$ & $12$ & $81$ & $2.01\times10^{-4}$ & $2.02\times10^{-4}$ & $114.05$ s & $1124\times$ \\
$[8.5,10]$ & $33$ & $3$  & $65$ & $1.68\times10^{-4}$ & $1.66\times10^{-4}$ & $45.63$ s & $1952\times$ \\
$[13,15]$  & $45$ & $15$ & $89$ & $2.41\times10^{-4}$ & $2.42\times10^{-4}$ & $138.47$ s & $879\times$ \\
\bottomrule
\end{tabular}
\caption{Performance of the ROM on three frequency bands of the Marmousi model.}
\label{tab:bands}
\end{table}

In the three cases the greedy loop reaches the prescribed tolerance, and the
maximum relative error on the test set remains below
$\varepsilon_{\mathrm{tol}}$. The two additional bands are wider than the
reference one ($2$ Hz against $1.5$ Hz) and require more greedy iterations,
$12$ and $15$ against $3$, and larger reduced bases, $n_U=n_\Lambda=81$ and
$89$ against $65$. The offline cost increases accordingly, from
$45.63\,\mathrm{s}$ to $114.05\,\mathrm{s}$ and $138.47\,\mathrm{s}$. The
online speedup decreases with the size of the reduced bases, from
$1952\times$ to $1124\times$ and $879\times$. These results also show the influence of the frequency range on ROM performance, we observe that since higher frequencies yield more oscillatory solutions, it require larger reduced bases to achieve the prescribed accuracy.

\section{Three-dimensional model}
\label{sec:3d}

We now assess the frequency-domain ROM construction of Section~\ref{sec:ROM} on a
three-dimensional heterogeneous medium. We consider the acoustic wave
propagation problem \eqref{eq:acoustic_strong} on a rectangular domain
$\Omega\subset\mathbb{R}^3$, with a Neumann boundary condition imposed on the
top boundary and absorbing boundary conditions on the remaining five faces.
The problem is discretized using the HDG scheme described in
Section~\ref{sec:discretization}, with polynomial degrees
$p=3$, on a tetrahedral mesh comprising
$N_e=148{,}422$ elements leading to $N_U=11,873,760$ and $N_\Lambda=3,059,460$. The velocity model is heterogeneous, with constant
density $\rho=1$. A sharp localized Gaussian source is located at $\boldsymbol{x}= (1500, 700, 200)$. 

\begin{figure}[ht!]
    \centering
    \subfloat[3D wave speed model of size 2.54$\times$1.44$\times$1.22 km$^3$ and 
    section in $y=$1.14 km.]
    {\includegraphics[scale=1]{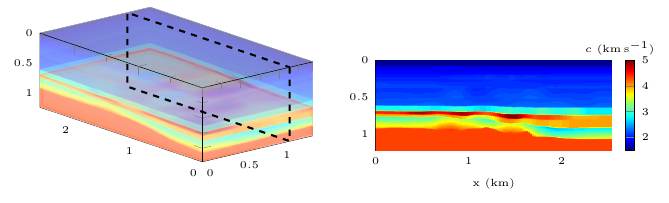}}

    \subfloat[
    Real part of the pressure field at 
    frequency $\omega=2\pi\times 8.02$
    (left)
    and $2\pi\times 10.92$
    (right)
    using the high-fidelity solver, for a source located at $(1.5,\,0.7,\,0.2)$~km.
    ]{
    \includegraphics[width=17em]{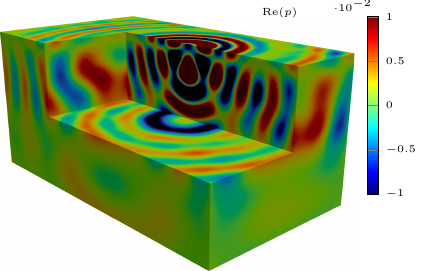} \hspace*{3em}
    \includegraphics[width=17em]{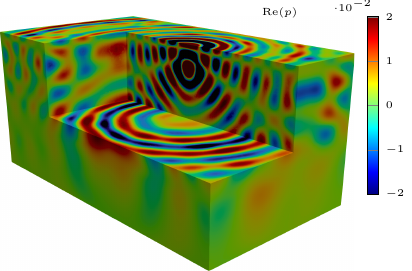}
    }

    \caption{Three-dimensional problem using a frequency interval from 8.02 to 10.92 Hz. }
    \label{fig:pressure_3D}
\end{figure}

The frequency band of interest is
$\mathcal{D}_\omega=[8,11]\times 2\pi$, with damping parameter $s=0.001$. For
simplicity of presentation, the angular frequency $\omega$ is identified with
its frequency in the remainder of this section. The reduced model is
constructed using Algorithm~\ref{alg:overall}. Figure~\ref{fig:pressure_3D}
illustrates the resulting high-fidelity pressure field $p$ at $\omega=8.02$ and $\omega=10.92$,
highlighting the heterogeneous nature of the wavefield.

The numerical investigation follows the same organization as in
Section~\ref{sec:marmousi}: we first examine separately the sensitivity
enrichment, the error bounds and the error indicator, and then assess the
convergence, cost and accuracy of the full greedy procedure. The error metrics
$E_{\mathrm{proj}}$, $E_{\mathrm{rel}}$ and $\theta_{\mathrm{eff}}$ are those
defined in \eqref{eq:proj_ern_U}--\eqref{eq:eff_index}.

\paragraph{Computational setup.}
All 3D computations were performed on PlaFRIM platform \cite{PlaFRIM}, on nodes with two 24-core AMD Zen4 processors and 256\,GB of memory. We use a hybrid MPI/OpenMP parallelization (Cf. \cref{sec:overall_algo}) with $4$ MPI processes per node and $12$ OpenMP threads per process. The studies of the enrichment, the gaps and the error indicator run on a single node ($4$ MPI processes), while the greedy runs and the final ROM use $4$ nodes ($16$ MPI processes).

\subsection{Enhancement of the ROB via frequency derivatives}

We assess the benefit of enriching the training snapshots with their first order
frequency derivatives, with $n_{\mathrm{train}}=20$ frequency points uniformly distributed over
$\mathcal{D}_\omega$, and an out-of-sample test set of $n_{\mathrm{test}}=15$
frequencies.

At the same POD tolerance ($10^{-10}$), enrichment increases the retained
basis size from $n_U=n_\Lambda=20$ to $n_U=n_\Lambda=35$, indicating that the
derivative information is not redundant with the primal snapshots and
meaningfully enriches the reduced space. Table~\ref{tab:enrichment_3d} reports
the relative projection error $E_{\mathrm{proj}}$, the relative $L^2$ error
$E_{\mathrm{rel}}$, and the effectivity index $\theta_{\mathrm{eff}}$ for both
variables, with and without enrichment.

We observe an improvement of nearly two orders of magnitude in the projection
error for both variables $u_h$ and $\lambda_h$: the mean projection error for
$u_h$ decreases from $1.30\times10^{-3}$ to $1.80\times10^{-5}$. A comparable
improvement is observed in the ROM error: the maximum error in $u_h$ decreases
from $9.05\times10^{-3}$ to $1.01\times10^{-4}$, with $\lambda_h$ exhibiting
the same behavior.

This improvement comes at a negligible additional offline cost. Computing the
$20$ high-fidelity snapshots together with their frequency derivatives
requires $2867\,\mathrm{s}$ in total, of which only $44.7\,\mathrm{s}$,
corresponding to approximately $1.6\%$ of the total training cost, is
associated with the derivative computation, the same LU factorization being
reused as described in \cref{sec:basis_construction}.

The effectivity index remains in the same range with and without enrichment,
with $\theta_{\mathrm{eff}}\approx1.5$--$1.8$ in both cases. The Galerkin
reduced solution thus remains quasi-optimal with respect to the
best-approximation error, which justifies the Galerkin formulation without
resorting to a Petrov--Galerkin projection. We therefore retain the
frequency-sensitivity enrichment for the remainder of the tests.

\begin{table}[ht!]
\centering
\small
\begin{tabular}{lcccccc}
\toprule  & rank $n$
 & mean $E_{\mathrm{proj}}$ & max $E_{\mathrm{proj}}$ & mean $E_{\mathrm{rel}}$ & max $E_{\mathrm{rel}}$ & $\theta_{\mathrm{eff}}$ (mean/max) \\
\midrule
\multicolumn{7}{l}{\textit{Variable $u_h$}} \\
without freq. sensitivity & $20$ & $1.30\times10^{-3}$ & $5.86\times10^{-3}$ & $2.03\times10^{-3}$ & $9.05\times10^{-3}$ & $1.59$ / $1.54$ \\
with freq. sensitivity    & $35$ & $1.80\times10^{-5}$ & $6.76\times10^{-5}$ & $3.08\times10^{-5}$ & $1.01\times10^{-4}$ & $1.80$ / $1.49$ \\
\addlinespace
\multicolumn{7}{l}{\textit{Variable $\lambda_h$}} \\
without freq. sensitivity & $20$ & $1.30\times10^{-3}$ & $5.88\times10^{-3}$ & $2.03\times10^{-3}$ & $9.06\times10^{-3}$ & $1.58$ / $1.54$ \\
with freq. sensitivity    & $35$ & $1.81\times10^{-5}$ & $6.78\times10^{-5}$ & $3.07\times10^{-5}$ & $1.01\times10^{-4}$ & $1.79$ / $1.49$ \\
\bottomrule
\end{tabular}
\caption{Enhancement of the approximation properties using frequency derivatives (3D heterogeneous test case, $n_{\mathrm{train}}=20$).}
\label{tab:enrichment_3d}
\end{table}

\subsection{Analysis of error bounds}

We repeat in 3D the study of the gap terms of Proposition~\ref{prop:trace},
following the same protocol as in Section~\ref{sec:marmousi}. The trace basis
$\boldsymbol{Z}_\Lambda$ is kept fixed throughout the experiment, with a POD
tolerance of $10^{-10}$, resulting in 
$E^\Lambda_{\rm proj}=1.1865\times10^{-5}$. The volume basis $\boldsymbol{Z}_U$
is varied by changing its POD tolerance ${\rm tol}_U$. For each tolerance, we compare the ROM obtained without lifting to the ROM enriched with the lifting condition at $n_{\rm lift}=3$ selected frequencies. The resulting gap terms, projection
errors, and relative errors are reported in Table~\ref{tab:gaps_3d}.

\begin{table}[H]
\centering
\small
\setlength{\tabcolsep}{4pt}
\begin{tabular}{lccccccc}
\toprule
 ${\rm tol}_U$ & $\delta_K$ & $\delta_F$
 & $E^{U}_{\mathrm{proj}}$ & $E^{\Lambda}_{\mathrm{proj}}$
 & $E^{U}_{\mathrm{rel}}$ & $E^{\Lambda}_{\mathrm{rel}}$
 & $\theta^{U}_{\mathrm{eff}}$ \\
\midrule
\multicolumn{8}{l}{\textit{without lifting}} \\
$10^{-6}$  & $1.00$              & $7.44\times10^{-3}$ & $2.9956\times10^{-4}$ & $1.1865\times10^{-5}$ & $9.7607\times10^{-4}$ & $9.7551\times10^{-4}$ & $3.3$  \\
$10^{-8}$  & $1.00$              & $5.14\times10^{-3}$ & $4.7048\times10^{-5}$ & $1.1865\times10^{-5}$ & $1.0968\times10^{-4}$ & $1.0963\times10^{-4}$ & $2.6$  \\
$10^{-10}$ & $2.05\times10^{-4}$ & $5.41\times10^{-5}$ & $1.1840\times10^{-5}$ & $1.1865\times10^{-5}$ & $2.4412\times10^{-5}$ & $2.4401\times10^{-5}$ & $2.1$  \\
$10^{-12}$ & $2.05\times10^{-4}$ & $5.48\times10^{-5}$ & $1.5209\times10^{-6}$ & $1.1865\times10^{-5}$ & $1.9968\times10^{-5}$ & $1.9954\times10^{-5}$ & $14.2$ \\
$10^{-14}$ & $2.05\times10^{-4}$ & $5.38\times10^{-5}$ & $2.9819\times10^{-7}$ & $1.1865\times10^{-5}$ & $1.9969\times10^{-5}$ & $1.9954\times10^{-5}$ & $93.8$ \\
$10^{-16}$ & $2.05\times10^{-4}$ & $5.32\times10^{-5}$ & $1.0064\times10^{-7}$ & $1.1865\times10^{-5}$ & $1.9969\times10^{-5}$ & $1.9955\times10^{-5}$ & $618.5$ \\
\addlinespace
\multicolumn{8}{l}{\textit{with lifting}, $n_{\rm lift}=3$} \\
$10^{-6}$  & $5.38\times10^{-12}$ & $4.86\times10^{-12}$ & $3.5826\times10^{-6}$ & $1.1865\times10^{-5}$ & $1.8970\times10^{-5}$ & $1.8952\times10^{-5}$ & $5.5$ \\
$10^{-8}$  & $3.53\times10^{-12}$ & $1.26\times10^{-12}$ & $1.6020\times10^{-6}$ & $1.1865\times10^{-5}$ & $1.8970\times10^{-5}$ & $1.8952\times10^{-5}$ & $12.2$ \\
$10^{-10}$ & $2.66\times10^{-12}$ & $5.32\times10^{-13}$ & $7.4174\times10^{-7}$ & $1.1865\times10^{-5}$ & $1.8970\times10^{-5}$ & $1.8952\times10^{-5}$ & $29.7$ \\
$10^{-12}$ & $1.99\times10^{-12}$ & $3.30\times10^{-13}$ & $3.9703\times10^{-7}$ & $1.1865\times10^{-5}$ & $1.8970\times10^{-5}$ & $1.8952\times10^{-5}$ & $59.7$ \\
$10^{-14}$ & $1.51\times10^{-12}$ & $2.88\times10^{-13}$ & $2.4362\times10^{-7}$ & $1.1865\times10^{-5}$ & $1.8970\times10^{-5}$ & $1.8952\times10^{-5}$ & $108.1$ \\
$10^{-16}$ & $1.47\times10^{-12}$ & $2.89\times10^{-13}$ & $8.6017\times10^{-8}$ & $1.1865\times10^{-5}$ & $1.8970\times10^{-5}$ & $1.8952\times10^{-5}$ & $822.2$ \\
\bottomrule
\end{tabular}
\caption{Gaps $\delta_K$, $\delta_F$, projection errors and relative errors, with and without lifting, at frozen $\boldsymbol{Z}_\Lambda$, 3D heterogeneous case.}
\label{tab:gaps_3d}
\end{table}

The lifted configurations confirm the decoupling predicted by
Proposition~\ref{prop:trace}, as observed on the Marmousi model. Enforcing the
lifting condition reduces both gap terms to numerical precision, with
$\delta_K\approx10^{-12}$ and $\delta_F\approx10^{-13}$--$10^{-12}$ throughout
the experiment. Consequently, the trace error becomes independent of the
volume reduced space: $E^\Lambda_{\rm rel}$ remains equal to
$1.8952\times10^{-5}$ for all six lifted configurations, whereas it varies
with the volume basis through the gap terms in the unlifted configurations.
Consistently with \cref{eq:sum_lambda}, the trace error is then governed only
by the approximation properties of $\boldsymbol{Z}_\Lambda$, giving
$\theta^\Lambda_{\rm eff}=1.58$.

The lifted configurations also illustrate the behavior predicted by
Proposition~\ref{prop:primal}. The relative volume error remains constant at
$E^U_{\rm rel}=1.8970\times10^{-5}$ across all volume-basis tolerances, while
the volume projection error is more than one order of magnitude smaller
($0.09$--$3.6\times10^{-6}$). The volume error is therefore controlled by the
trace error propagated to the volume variable, in agreement with the structure
of the bound~\eqref{eq:sum_primal} which explains the large values of
$\theta^U_{\rm eff}$.  

The unlifted configurations reproduce the observations made in 2D. As the accuracy of the reduced basis $\boldsymbol{Z}_U$ increases, the gap terms $\delta_K$ and $\delta_F$ decrease yielding comparable errors between the lifted and the unlifted configurations. As discussed in Section~\ref{sec:marmousi}, this
follows from the structure of the volume snapshots, which already contain the
lifting of the corresponding trace snapshots. 
We therefore retain Algorithm~\ref{alg:overall} as it is, without explicitly
enforcing the lifting condition, as in the 2D case.

\subsection{Behavior of the error indicator}

\begin{figure}[H]
\centering
\includegraphics[scale=0.35]{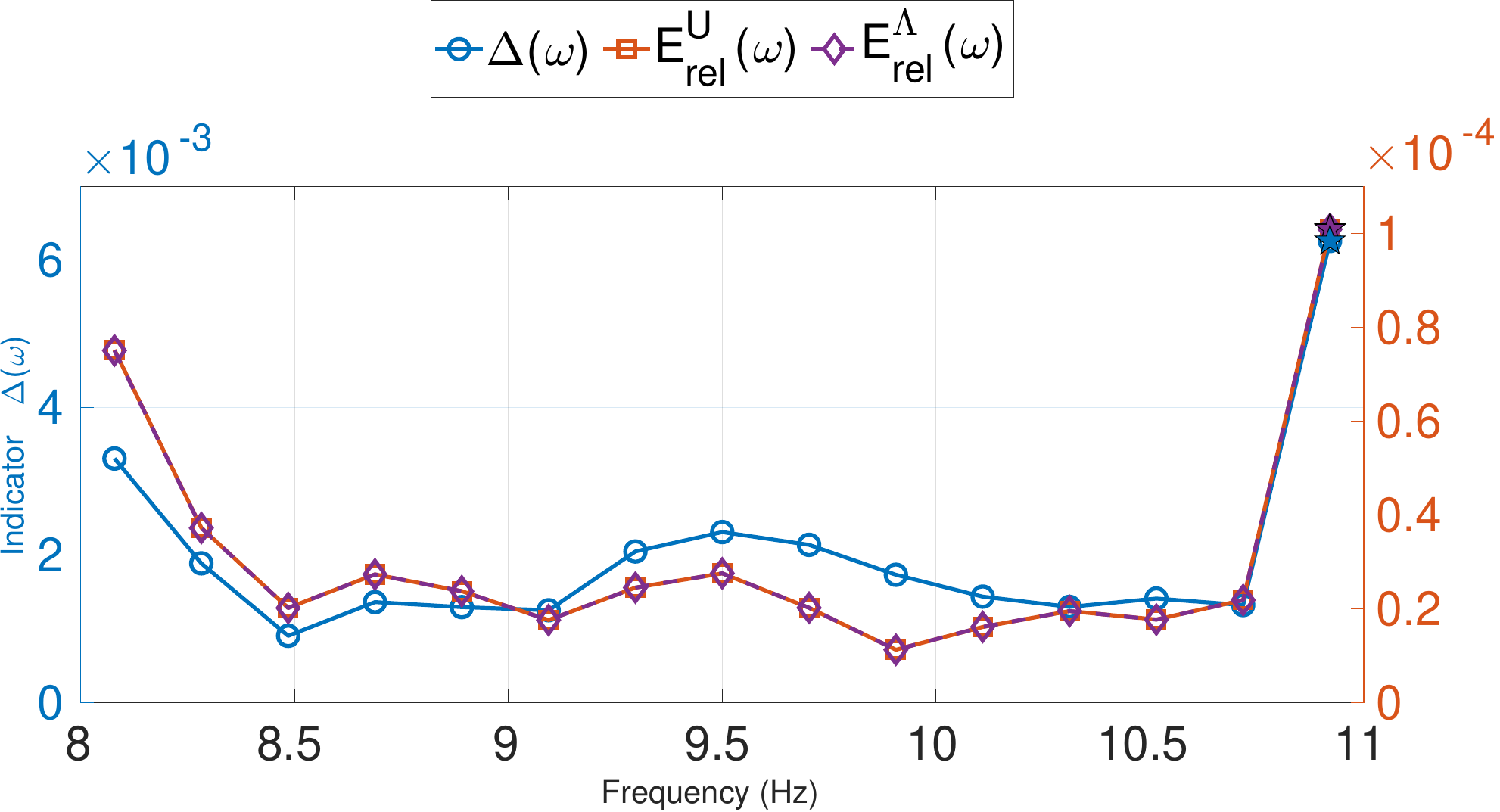}
\caption{Indicator $\Delta(\omega)$ and true relative error
$E^U_{\mathrm{rel}}(\omega)$ on the test set; the marked points indicate $\arg\max_\omega
\Delta(\omega)$ and $\arg\max_\omega E^U_{\mathrm{rel}}(\omega)$.}
\label{fig:indicator_test1_3d}
\end{figure}

We evaluate the error indicator $\Delta(\omega)$ defined in
\eqref{eq:greedy_indicator}. On a plain POD run with a uniform training grid
of $n_{\mathrm{train}}=20$ frequencies over $\mathcal{D}_\omega$, we compare
$\Delta(\omega)$ with the relative error $E_{\mathrm{rel}}(\omega)$ defined in
\cref{eq:rel_ern_U}--\cref{eq:rel_err_L}, on a test set of
$n_{\mathrm{test}}=15$ frequencies. Figure~\ref{fig:indicator_test1_3d} shows
that the residual-based indicator follows the behavior of the relative error
over the whole band and reproduces its local variations.

Both quantities attain their maximum at the same frequency, $\omega=10.92$,
so that the frequency selected by the indicator in the greedy step is the most informative one.

As in 2D, the relative errors of $u_h$ and $\lambda_h$ coincide, since the
skeleton variable $\lambda_h$ represents the numerical trace of the pressure
component of $u_h$. The indicator is again to be interpreted as a tool for
locating the frequencies of largest ROM error rather than as a sharp estimate
of its magnitude (see Section~\ref{sec:marmousi}).

\subsection{Convergence of the greedy algorithm}

We now assess the overall performance of \cref{alg:overall}. Starting from
$n_{\mathrm{init}}=15$ uniformly distributed snapshots and a candidate set of
$n_{\mathrm{cand}}=100$ frequencies, the greedy loop is terminated when the
relative error falls below $\varepsilon_{\mathrm{tol}}=10^{-4}$.

We report the evolution of the indicator $\Delta(\omega)$ and of the true
relative error $E_{\mathrm{rel}}(\omega)$ at each iteration, to assess whether
reducing the indicator along the greedy loop effectively reduces the true
error. The loop terminates after five iterations.

\begin{figure}[H]
\centering
\includegraphics[scale=0.35]{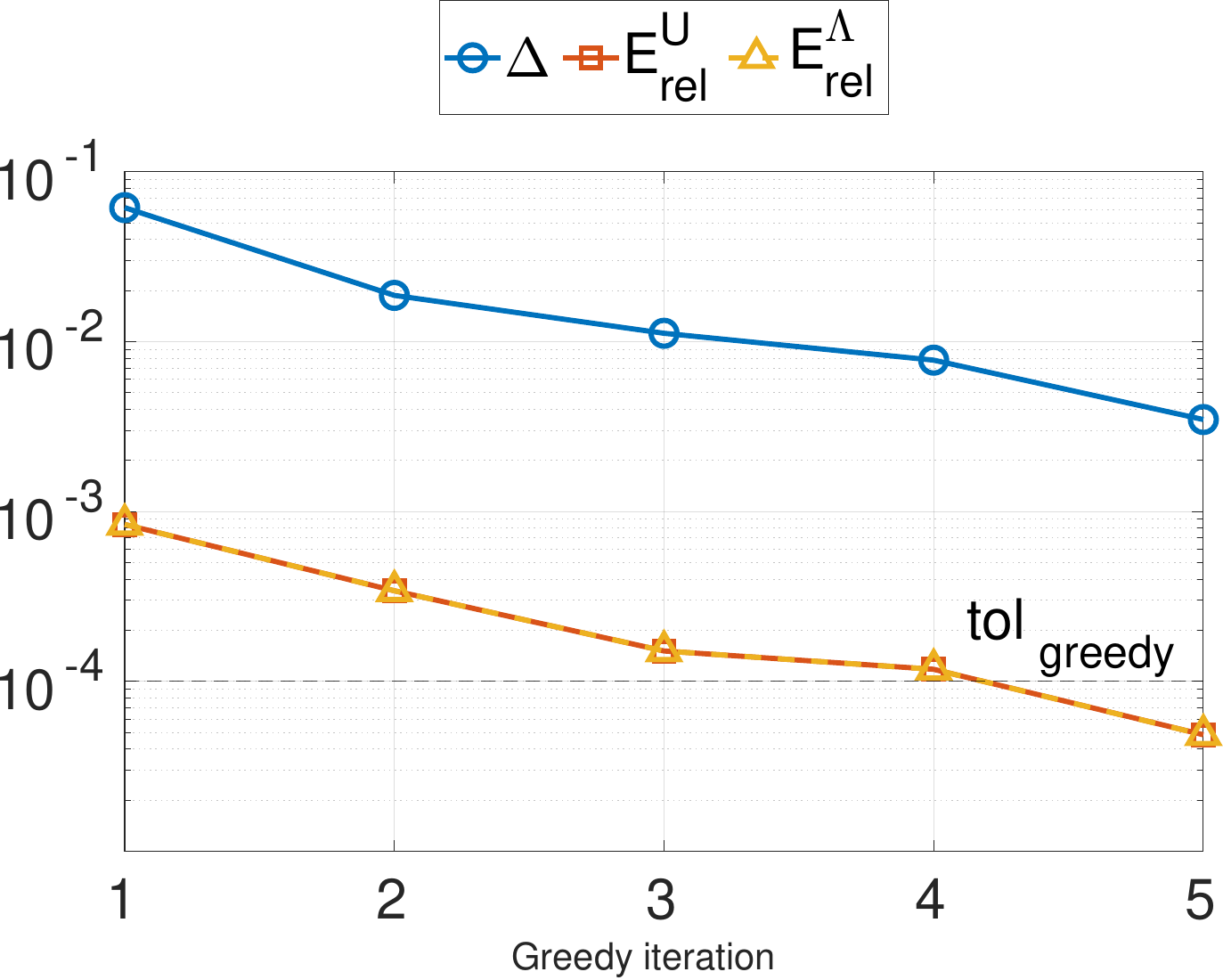}
\caption{Evolution of the indicator $\Delta$ and the relative errors
$E^U_{\mathrm{rel}}$, $E^\Lambda_{\mathrm{rel}}$ along the greedy
iterations, the loop stops when $E^U_{\rm rel}\, , \,  E^\Lambda_{\rm rel} \le \varepsilon_{\mathrm{tol}}= 10^{-4}$.}
\label{fig:indicator_test2_3d}
\end{figure}

Figure~\ref{fig:indicator_test2_3d} shows that both $\Delta$ and
$E_{\mathrm{rel}}$ decrease monotonically, with the relative error decreasing
from $8.40\times10^{-4}$ to $4.85\times10^{-5}$. Evaluating the indicator over
the $100$ candidates costs $1.1\times10^{-3}\,\mathrm{s}$ per iteration, a
negligible fraction of the iteration cost, which makes a dense exploration of
the parameter domain affordable.

\subsection{Offline and online costs}

Offline, constructing the ROM requires $20$ high-fidelity (HF) solves
($n_{\mathrm{init}}=15$ initial snapshots plus $5$ greedy iterations).
Together with their frequency derivatives, these yield $40$ candidate vectors,
compressed by POD into a reduced basis of size $n_U=n_\Lambda=39$. The HF
solves account for $1070.20\,\mathrm{s}$ ($808\,\mathrm{s}$ for the initial
training set and $262\,\mathrm{s}$ within the loop), out of a total offline
cost of $1289.17\,\mathrm{s}$. The remaining $218.97\,\mathrm{s}$ is spent on
POD compression, Gram matrix assembly, projection, and indicator evaluation
across the greedy loop (cf.\ \cref{sec:basis_construction}). The offline cost
is therefore dominated by the HF solves, while the compression and projection
procedures account for only about $17\%$ of the total offline cost.

Online, evaluating the ROM at a new frequency costs on average
$1.11\times10^{-4}\,\mathrm{s}$, compared with $48.25\,\mathrm{s}$ for a
full-order solve, yielding a speedup of approximately $4.3\times10^{5}$
(Table~\ref{tab:speedup_3d}). The offline training is therefore recovered
after approximately $27$ online queries.

\begin{table}[htbp]
\centering
\begin{tabular}{lc}
\toprule
Offline: number of HF solves & $20$ \\
Offline: HF solve time & $1070.20$ s \\
Offline: POD + projection + residual indicator build & $218.97$ s \\
Offline: total time & $1289.17$ s \\
\midrule
Online: FOM time per query & $48.25$ s \\
Online: ROM time per query & $1.11\times10^{-4}$ s \\
Online: speedup per query & $4.3\times10^{5}\times$ \\
\bottomrule
\end{tabular}
\caption{Offline and online computational costs of the final greedy ROM ($5$ iterations, $n_U=n_\Lambda=39$), evaluated on $n_{\mathrm{test}}=50$ independent test frequencies.}
\label{tab:speedup_3d}
\end{table}

\subsection{Accuracy: FOM vs ROM}

The final ROM, obtained after $5$ greedy iterations with $n_U=n_\Lambda=39$, is
validated on an independent test grid of $n_{\mathrm{test}}=50$ frequencies
uniformly sampled over $\mathcal{D}_\omega$, disjoint from the training
frequency set. Table~\ref{tab:accuracy_3d} reports the projection errors, with
a mean value of $6.6\times10^{-6}$, indicating that the greedy procedure
selects snapshots that yield an accurate reduced basis for this oscillatory
and heterogeneous three-dimensional configuration.

Table~\ref{tab:accuracy_3d} also reports the relative $L^2$ errors over the
test set, separately for $(p_h,\boldsymbol{v}_h)$ and $\lambda_h$. For $p_h$,
the mean (resp. maximum) relative error is $1.04\times10^{-5}$
(resp. $4.03\times10^{-5}$), and the velocity components $v_{h,x}$, $v_{h,y}$ and
$v_{h,z}$ have errors of the same order of magnitude. The corresponding values
for $\lambda_h$ are $1.04\times10^{-5}$ ($4.02\times10^{-5}$). In all cases,
the maximum relative errors remain below the greedy tolerance
$\varepsilon_{\mathrm{tol}}=10^{-4}$ used to terminate the training. Finally,
the ratio between the relative and projection errors remains close to $1.6$
for all fields, confirming the quasi-optimality of the Galerkin ROM with
respect to the reduced bases.

\begin{table}[H]
\centering
\begin{tabular}{lcccc}
\toprule
 & mean $E_{\mathrm{proj}}$ & max $E_{\mathrm{proj}}$ & mean $E_{\mathrm{rel}}$ & max $E_{\mathrm{rel}}$ \\
\midrule
$p_h$      & $6.56\times10^{-6}$ & $2.51\times10^{-5}$ & $1.04\times10^{-5}$ & $4.03\times10^{-5}$ \\
$v_{h,x}$  & $6.34\times10^{-6}$ & $2.46\times10^{-5}$ & $9.87\times10^{-6}$ & $3.91\times10^{-5}$ \\
$v_{h,y}$  & $4.75\times10^{-6}$ & $1.84\times10^{-5}$ & $7.78\times10^{-6}$ & $3.12\times10^{-5}$ \\
$v_{h,z}$  & $8.05\times10^{-6}$ & $3.11\times10^{-5}$ & $1.08\times10^{-5}$ & $4.20\times10^{-5}$ \\
$\lambda_h$  & $6.58\times10^{-6}$ & $2.52\times10^{-5}$ & $1.04\times10^{-5}$ & $4.02\times10^{-5}$ \\
\bottomrule
\end{tabular}
\caption{Projection and relative $L^2$ errors of the final greedy ROM on an independent test grid of $n_{\mathrm{test}}=50$ frequencies over $\mathcal{D}_\omega=[8,11]\times 2\pi$.}
\label{tab:accuracy_3d}
\end{table}

\begin{figure}[ht!] \centering
  \includegraphics[width=15em]{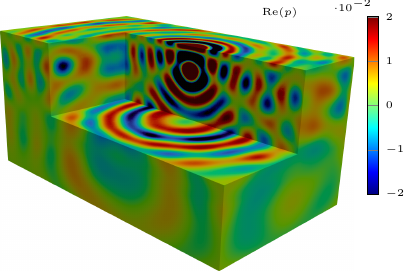} \hspace{1em}
  \includegraphics[width=15em]{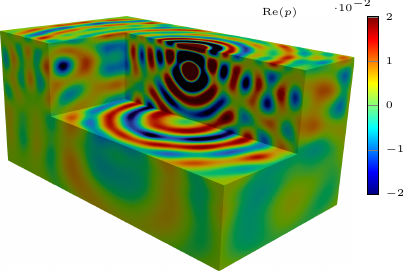}\hspace{1em}
  \includegraphics[width=15em]{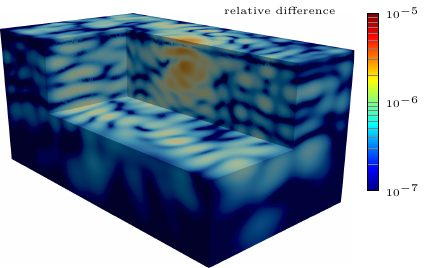}
  \caption{Real part of the pressure field for the high-fidelity solution (left) and
  the ROM prediction (middle) at the worst test frequency $\omega=2\pi\times 9.86$ and the relative
  difference (right).}
  \label{fig:3d-comparison-rom}
\end{figure}

Figure~\ref{fig:3d-comparison-rom} compares the high-fidelity and reduced pressure field at the test frequency of largest error, $\omega=9.86\times2\pi$. We note that the fields are visually identical. The associated relative error ranges from  $10^{-7}$ to $10^{-5}$, it is largest in the vicinity of the source and remains relatively small throughout the rest of the computational domain.

\subsection{Influence of the frequency band}

We repeat the greedy construction on two other frequency bands,
$\mathcal{D}_\omega=[4,7]\times2\pi$ and $\mathcal{D}_\omega=[13,15]\times2\pi$,
with the same setting as the previous band. The band $[13,15]$ is highly oscillatory, we  thus use polynomial degree $p=4$, leading to $N_U\approx20.8\times10^{6}$ and
$N_\Lambda\approx4.6\times10^{6}$. We note thus that the CPU costs are not directly
comparable with those of the two other bands since we do not use the same polynomial order, the experiments aims to assess how the reduced model depends on the frequency band. Table~\ref{tab:bands_3d} reports, for each band, the number of HF solves, the
size of the reduced bases, the maximum relative errors over the test set, the
offline cost and the online speedup.

\begin{table}[H]
\centering
\begin{tabular}{lcccccccc}
\toprule
band (Hz) & $p$ & HF solves & iterations & $n_U=n_\Lambda$ & max $E^U_{\mathrm{rel}}$ & max $E^\Lambda_{\mathrm{rel}}$ & offline time & speedup \\
\midrule
$[4,7]$   & $3$ & $18$ & $3$ & $34$ & $3.13\times10^{-5}$ & $3.13\times10^{-5}$ & $959$ s & $5.4\times10^{5}\times$ \\
$[8,11]$  & $3$ & $20$ & $5$ & $39$ & $4.20\times10^{-5}$ & $4.02\times10^{-5}$ & $1289$ s & $4.3\times10^{5}\times$ \\
$[13,15]$ & $4$ & $16$ & $1$ & $31$ & $3.49\times10^{-5}$ & $3.29\times10^{-5}$ & $2401$ s & $1.7\times10^{6}\times$ \\
\bottomrule
\end{tabular}
\caption{Greedy ROM performance on three frequency bands of the 3D heterogeneous model}
\label{tab:bands_3d}
\end{table}

In the three cases the greedy loop reaches the prescribed tolerance with at
most five iterations, and the maximum relative error on the test set remains
below $\varepsilon_{\mathrm{tol}}=10^{-4}$, with reduced bases of size
$n_U=n_\Lambda=31$ to $39$. The offline cost is dominated by the HF solves,
which account for $83\%$ to $97\%$ of the total. The cost of a reduced solve
is of the order of $10^{-4}\,\mathrm{s}$ in all cases , so that the speedup increases with the
cost of the full-order solve, from about $5\times10^{5}$ for $p=3$ to
$1.7\times10^{6}$ for $p=4$. The offline training is recovered after
approximately $20$, $27$ and $17$ online queries, respectively.

\section{Conclusion}
\label{sec:conclusion}

In this work, we have proposed a methodology for constructing reduced-order
models tailored to the structure of HDG discretizations. The methodology
combines several tools from the literature within a unified framework:
greedy sampling of the frequency band driven by a residual-based error
indicator, enrichment of the snapshots with frequency derivatives, and
POD compression performed separately for the volume and trace unknowns.

For the damped time-harmonic acoustic problem, we have established the
well-posedness of both the HDG formulation and its reduced counterpart,
with the same stability constant. This result relies on constructing
decoupled reduced bases for the volume and trace unknowns. The resulting
stability constant is explicit and characterizes the dependence of the
stability on the damping, the medium parameters, and the mesh. We have
also investigated different projection strategies and their impact on the
reduced model. In particular, we have analyzed and quantified the
interaction between projection and static condensation, which do not
commute in general. This analysis has led to error estimates and the
identification of a lifting condition that eliminates the resulting gap.

The error estimates have been assessed numerically in two and three
dimensions. The results show that lifting decouples the trace error from
the approximation error associated with the volume reduced basis, thereby
recovering the structure of the HDG discretization. Moreover, we show that the
volume error is bounded in terms of the propagated trace error regardless of the accuracy of its reduced basis. 
As the volume basis becomes sufficiently accurate, however, the resulting
errors become comparable to those obtained with lifting. For the
problems considered here, lifting is therefore not necessary in practice,
provided that the volume reduced basis is sufficiently rich.

Finally, the methodology has been implemented in the open-source \texttt{Hawen}
code \cite{faucher2021hawen}, using a hybrid MPI/OpenMP parallelization
that exploits the element-wise structure of the HDG scheme to make
high-dimensional problems tractable. The methodology has been validated on realistic
two- and three-dimensional complex heterogeneous media. For the Marmousi model,
a reduced model with $65$ modes per variable achieves relative errors
below $5\times10^{-4}$ and an online speedup of approximately
$2\times10^{3}$. In this case, the offline cost is recovered after $71$ online queries. In three dimensions, for a problem with approximately
$12$ million volume unknowns, only $39$ modes are required to achieve
relative errors below $10^{-4}$, with an online speedup of approximately
$4\times10^{5}$. In this case, the offline cost is recovered after fewer
than $30$ online queries.

Future work will investigate the extension of the proposed methodology to elastic wave propagation and coupled acoustic–elastic problems, where the interface conditions coupling the acoustic and elastic fields must also be incorporated into the reduced model. We further consider extensions to full waveform inversion (FWI), which raises two challenges for model reduction. First, the reduced space is currently constructed for a single source, whereas FWI typically involves multiple sources. Constructing reduced spaces that remains accurate with multiple source configurations is therefore an important direction for future work. Second, the medium parameters change throughout the inversion process, requiring the reduced basis to account for medium variability while preserving the accuracy of the reduced-order model. This motivates the development of efficient strategies to adapt the reduced basis as the medium evolves, without rebuilding the reduced-order model at each iteration.

\section{Acknowledgment}

This work was supported by the EXAMA (Methods and Algorithms 
at Exascale) project under grant ANR-22-EXNU-0002.
F.F. and I.T. acknowledge funding by the European Union
with ERC Project \textsc{Incorwave} -- grant 101116288.
Views and opinions expressed are however those of the authors
only and do not necessarily reflect those
of the European Union or the European Research Council
Executive Agency (ERCEA). Neither the European Union nor the
granting authority can be held responsible for them.

\bibliographystyle{abbrv}
\bibliography{references}
\end{document}